%% file: huang_and_hardle_jbes.tex
\documentclass[12pt, a4paper]{article}

\usepackage{amsmath,rotating,graphicx,latexsym,verbatim,amssymb}
\usepackage{setspace}
\usepackage{times}
\usepackage{bm}
\usepackage{natbib}
\usepackage{amsthm}
\usepackage{caption}
\usepackage[plain,noend]{algorithm2e}

\input{defs}
\renewcommand{\baselinestretch}{1.45}
\newtheorem{theorem}{Theorem}

\def\singlespace{\def\baselinestretch{1}\@normalsize}

\newcommand{\by}{\mbox{\bf y}}

\newcommand{\bb}{\mbox{\bf b}}

\newcommand{\bz}{\mbox{\bf z}}

\newcommand{\bx}{\mbox{\bf x}}

\newcommand{\ba}{\mbox{\bf a}}

\newcommand{\bbeta}{\boldsymbol\beta}
\newcommand{\balpha}{\boldsymbol\alpha}
\newcommand{\boldmu}{{\mbox{\boldmath $\mu$}}}
\newcommand{\btheta}{{\mbox{\boldmath $\theta$}}}

\newtheorem*{proposition}{Proposition}
\begin{document}

\setcounter{page}{1}
\renewcommand{\baselinestretch}{1.2}

\begin{center}
{\Large \bf  Analysis of Deviance for Hypothesis Testing in Generalized Partially Linear Models  
  \vspace{0.3cm}\\}
 Wolfgang Karl H\"ardle\\
Center for Applied Statistics and Economics, Humboldt University, 10099 Berlin, GERMANY \\
{\tt haerdle@wiwi.hu-berlin.de}   \\
  AND\\      Li-Shan Huang\\
 Institute of Statistics, National Tsing Hua University, 30013,  TAIWAN \\
{\tt lhuang@stat.nthu.edu.tw} \\
\noindent  \today 
\end{center}

\begin{abstract}
In this study, we develop nonparametric analysis of deviance tools for generalized partially linear models based on local polynomial fitting.  Assuming a canonical link, we propose expressions for  both local  and global analysis of deviance, which admit an additivity property that reduces to  analysis of variance decompositions in the Gaussian case.  Chi-square tests based on integrated  likelihood functions  are proposed to formally  test whether the nonparametric term is significant. Simulation results are shown  to illustrate the proposed chi-square tests and to compare them with an existing procedure based on penalized splines. The methodology is applied to  German Bundesbank Federal Reserve data.
\end{abstract}


\section{Introduction}

Generalized linear models  (McCullagh and Nelder 1989) are a large class of statistical models for relating  a response variable to linear combinations of predictor variables. The models allow  the response variable to  follow probability distributions in the exponential family  such as the  Binomial and Poisson,  generalizing  the Gaussian distribution in linear models, though
a major limitation  is   the
prespecified  linear form of predictors. 
Generalized partially linear models   (Green and Silverman1994; Carroll et al. 1997; H\"ardle et al. 2004) allow for a nonparametric component for a continuous covariate while retaining the ease of linear relationships for the remaining variables. It is more flexible than the conventional linear approach
 and is a special case of  generalized additive models  (Hastie and Tibshirani 1990; Wood 2006) which include multiple nonparametric components. 
H\"ardle et al. (1998)  applied  the generalized partially linear model to 1991 East-West German migration data  to model the probability of  migration with a nonlinear relationship to household income and linear relationships to other covariates such as age, gender, and employment status. Wood (2006, p. 248) gave an example of modeling the daily  total deaths in Chicago  in the period 1987-2000 as a Poisson distribution with a nonlinear trend of time and linear  effects of daily temperature and daily air-pollution levels of ozone, sulfur dioxide, and pm10.
An illustrating finance example in Section 6 of this paper is on bankruptcy prediction  for firms, known as rating or scoring,   from a set of financial ratio variables. The logistic partially linear model is used to model the  probability of default  with a  nonlinear relationship to the account payable turnover  ratio, which is  a short-term liquidity measure, and linear relationships to some selected financial ratios.

 In applying generalized partially linear models to data, inference tools to examine  whether the nonparametric term is significant  are of interest. For example, in H\"ardle et al. (1998), the nonlinear estimated  function of household income showed a saturation in the intention to migrate for higher income households and the question was whether 
the overall income  effect was significant statistically. 
As analysis of deviance was developed for generalized linear models  (McCullagh and Nelder 1989), it is natural to ask whether one can extend it for  generalized partially linear models.
Though Hastie and Tibshirani (1990) briefly  discussed analysis of deviance for generalized additive models, they noted that ``the distribution theory, however, is undeveloped"  and ``informal deviance tests with some heuristic justification" were adopted. The present paper fills the gap by establishing local and global  analysis of deviance expressions 
for  generalized partially linear models and  developing associated tests for checking whether the nonparametric term is significant.
 Li and Liang (2008) addressed assessing the significance of the nonparametric term  in the local polynomial settings  by extending generalized likelihood ratio tests (Fan et al. 2001), which have  asymptotic   chi-square distributions.   Wood (2013) discusses approximate $p$-values for testing significance of smooth components of semiparametric generalized additive models by Wald-type tests based on penalized splines. We remark that testing in the generalized partially linear models is relatively less developed as compared to the special case  of partially linear models  under  the Gaussian distribution (H\"ardle et al. 2004). Hence, there is a need for developing analysis of deviance  tools parallel to those in generalized linear models for applications of generalized partially linear models. 

Based on  the local polynomial  approach (Fan and Gijbels 1996) and assuming a canonical link in generalized partially linear models, we propose   local and global expressions for analysis of deviance, with the latter  obtained by integrating the corresponding local likelihood quantities.
This mimics the ``integrated likelihood" approach discussed by   Lehmann (2006) and Severini (2007). Though the idea of local likelihood  has been around for some time (Hastie and Tibshirani 1987; Loader 1999),  we are not aware of using the integrated likelihood approach to combine the information of local likelihood in the smoothing literature.   Then integrated likelihood ratio tests  with asymptotic chi-square distributions are proposed  to check whether the nonparametric term is significant. Our work extends the classic analysis of deviance  to generalized partially linear models  with theoretical justifications,  and generalizes the work of Huang and Chen (2008) and Huang and Davidson (2010) in  a special case of the Gaussian distribution to distributions in the canonical exponential family.

The organization of this paper is as follows.  Section 2 outlines the analysis of deviance for nested hypotheses in parametric generalized  linear models by Simon (1973). In Section 3, we propose local and global analysis of deviance for nonparametric  models   in Theorem 1 for the simpler case with the nonparametric term as the only predictor. By combining local likelihood through integration,  as a  by-product,  new estimators for the canonical parameter and response mean  are given in equation (12),  and Theorem 2 shows that the integrated likelihood quantities are asymptotically global likelihood quantities with the new estimators. Theorem 3 proposes  integrated likelihood ratio tests with asymptotic chi-square distributions for testing whether the nonparametric term is significant.  Section 4 presents an extension of Theorems 1-3 to  generalized partially linear models as Theorems 4 and 5.   In  Section 5, we  illustrate the potential  usefulness of the new tests with  simulated data   and compare with  tests  by Wood (2013)   in  the R package {\tt mgcv}.  Section 6 applies the methodology to  2002 German Bundesbank Federal Reserve data and Section 7 gives some concluding remarks and directions for future research.

\section{Preliminaries}

We first describe  generalized linear models  based on McCullagh and Nelder (1989). Let $(x_1, y_1), \dots, (x_n, y_n)$ be independent data pairs with the conditional density of $Y$ given covariate $X=x$ from a one-parameter exponential family:
\begin{equation}
L(y; \theta(x))=\exp\left[ \frac{y \theta(x)- b\{\theta(x)\}} { a(\phi)} + c(y,\phi) \right],
\label{glm}
\end{equation}
where  $a(\cdot)>0$, $b(\cdot)$ and $c(\cdot)$ are known functions,  $\phi$ is known or a nuisance parameter, and  $\theta$ is the canonical parameter with
the conditional mean of response,  $ \E(Y \mid X=x) =\mu(x)=b^{\prime}\{\theta(x)\}$.  
 A transformation of mean $G\{\mu(x)\}$ may be modelled linearly  by
$G\{\mu(x)\}= b_0 + b_1 x$, where  $G(\cdot)$ is called the ``link" function and estimates of $b_0$ and $b_1$ are obtained by maximum likelihood.
If $G(\cdot)= (b^{\prime})^{-1}(\cdot)$, then $G$ is the canonical link function that links $\theta$ to the linear predictor. For simplicity,  $G$ is the canonical link function throughout the paper and the dependence of  $\theta$ on covariates is often suppressed if no ambiguities result.

 Let $\ell(y; \theta)= \log L(y; \theta)$, $\hat{\theta}=G(\hat{\mu})$ denote the  fitted value of $\theta$ with corresponding $\hat{\mu}$, and $\tilde{\theta}=G(y)$ when the fitted value equals  the observed  $y$. The  deviance  $D$ (McCullagh and Nelder 1989), measuring the discrepancy between data $\by= (y_1, \dots, y_n)^{\top}$ and fitted $\hat{\boldmu}=(\hat{\mu}_1, \dots, \hat{\mu}_n)^{\top}$, 
is 
\begin{equation}
D(\by, \hat{\boldmu})= 2 \sum_i \{ y_i (\tilde{\theta}_i - \hat{\theta}_i) - b(\tilde{\theta}_i) + b(\hat{\theta}_i) \}.
\label{paradev}
\end{equation}
In the Gaussian case, $G$ is the identity link and $D=\sum_i (y_i -\hat{\mu}_i)^2$, which is
 the  residual sum of squares  in linear models. Let us now turn to  testing hypotheses about $\btheta=(\theta_1, \dots, \theta_n)^{\top}$.
 Assume that $D_j= \inf_{\btheta \in \Omega_j} D$,  $j=1,2$,  with $\Omega_2 \subseteq \Omega_1$.  The analysis of deviance usually refers to comparing two nested parametric models and  inference may be based on the difference $D_2- D_1$, which is simply the log likelihood ratio statistic with an asymptotic $\chi^2$ distribution. 
The conventional  analysis of deviance  is generally  not  parallel to the analysis of variance in linear models, in the sense that the former does not have all the  sum-of-squares quantities.

An attempt to mimic analysis of variance for (\ref{glm}) can be based on the  Kullback-Leibler (KL) 
divergence of two probability distributions with means $\boldmu_1$ and $\boldmu_2$:
\[
 KL(\boldmu_1, \boldmu_2) = 2 \E_{\boldmu_1}  \left[ \ell\{\by; G(\boldmu_1)\} -
 \ell \{\by; G(\boldmu_2)\} \right],
\]
where  $\boldmu_1$ and $\boldmu_2$ are  treated as fixed  values
and $\E_{\boldmu_1}$ is the conditional expectation with respect to $\by$  with $\boldmu=\boldmu_1$.
Simon (1973) showed that for nested  hypotheses $\Omega_2 \subset \Omega_1 \subset {\mathbb  R}^n$ with ${\mathbb  R}^n$ corresponding to the parameter space for an exact fit of $\tilde{\theta}$ and $\theta$  parameterized linearly in $\Omega_1$ and $\Omega_2$,  
\begin{equation}
 KL(\by, \boldmu_2)=  KL(\by, \boldmu_1) +  KL(\boldmu_1, \boldmu_2)
\label{simon}
\end{equation}
when plugging in the values of  maximum likelihood   estimates for  $\boldmu_1$ and $\boldmu_2$.  
In other words, (\ref{simon}) shows that the KL divergence exhibits the Pythagorean property. For the Gaussian distribution, 
 (\ref{simon})  reduces to the analysis of variance decomposition in linear models when  $\boldmu_1$ and $\boldmu_2$ correspond to  the linear   fit   and  the intercept-only model respectively, and the terms in (\ref{simon}) becomes total,  residual, and regression sums of squares   respectively.

A linear  form of $x$  may be restrictive and one may consider a nonparametric approach:  
\begin{equation}
G\{\mu(x)\}= m(x).
\label{nonpar}
\end{equation} 
Fan et al. (1995) discussed estimating  $m(\cdot)$ by maximizing a locally weighted likelihood with a local polynomial approximation. Based on  Taylor's expansion at $x$,
$\theta_i \approx \beta_0 + \beta_1 (x_i-x) + \dots + \beta_p (x_i -x)^p  \equiv 
 \theta_i(x)$. This approximation is plugged in the locally weighted  log-likelihood at $x$,  
\begin{equation}
\ell_x(\by; \btheta_x) \equiv \sum_i \ell\{y_i; \theta_i(x)\} K_h (x_i -x),
\label{lx}
\end{equation} 
where $\btheta_x=(\theta_1(x), \dots, \theta_n(x))^{\top},$  $K(\cdot)$ is usually  a   density function being symmetric at $0$,  $h$ is the bandwidth determining the neighborhood size, and $K_h(\cdot)= K(\cdot/h)/h$.   Then $\hat{\beta}=(\hat{\beta}_0, \dots, \hat{\beta}_p)^{\top}$  maximizing $\ell_x(\by; \btheta_x)$ is solved  and $j! \hat{\beta}_j$ estimates $m^{(j)}(x)$, $j=0, \dots, p$,  which is $\theta^{(j)}(x)$ with the canonical link.   Fan et al. (1995)   derived  asymptotic properties of $\hat{\beta}_j(x)$'s and adopted   $G^{-1}\{ \hat{\beta}_0(x)\}$ as an estimate for $\mu(x)$.  A further extension to
generalized partially linear models is 
\begin{equation}
G\{\mu(\bz, x)\}= \bz^{\top} \alpha + m(x),
\label{model}
\end{equation}
 where  $\bz$ is  a $K$-dimensional covariate vector. Without loss of generality, the intercept in (\ref{model}) is embedded in $m(\cdot)$. When $\alpha$ is unknown, estimation of $\alpha$ can be done via a two-step maximum likelihood procedure   that updates the linear and nonparametric estimates iteratively, as discussed in  Carroll et al. (1997), p. 479.

\section{Nonparametric Analysis of Deviance}

This section focuses on (\ref{nonpar}). 
We  start by deriving a local  analysis of deviance expression  for model (\ref{nonpar}) in the following by adapting (\ref{simon}) for locally weighted likelihood.
 Let $\hat{\theta}_i(x) = 
\hat{\beta}_0 + \dots + \hat{\beta}_p (x_i -x)^p$,  the resuting local polynomial  estimate of  $\theta_i$   at $x$, $\hat{\btheta}_x=(\hat{\theta}_1(x), \dots, \hat{\theta}_n(x))^{\top}$,  $\hat{\mu}_x(x_i)= G^{-1} \{\hat{\theta}_i(x) \}$, and $\hat{\boldmu}_x=(\hat{\mu}_x(x_1), \dots, \hat{\mu}_x(x_n))^{\top}$.
 As the  $\hat{\beta}_j$'s maximize $\ell_x(\by; \btheta_x)$,  the following equations hold:
\begin{eqnarray}
\nonumber \sum_i y_i (x_i -x)^j  K_h (x_i -x) & = & \sum_i \hat{\mu}_x(x_i)(x_i -x)^jK_h (x_i -x), \hspace{1cm} j = 0, \dots, p,\\
\sum_i y_i \hat{\theta}_i(x) K_h (x_i -x) & = & \sum_i \hat{\mu}_x(x_i) \hat{\theta}_i(x) K_h (x_i -x).
\label{eq1}
\end{eqnarray}
The last equation indicates that  $(\by - \hat{\boldmu}_x)$ is orthogonal to 
$\hat{\btheta}_x$ in the locally weighted inner product space with weights 
$K_h (x_i -x)$. Hence the fact of residuals being orthogonal to fitted values in ordinary linear models now becomes the fact of  local residuals $(\by - \hat{\boldmu}_x)$ being orthogonal to locally fitted canonical parameters $\hat{\btheta}_x$ in a kernel-weighted space.
For  $\ell_x(\by; \hat{\btheta}_x)$, an expression mimicking (\ref{paradev}) for  local deviance at $x$ is therefore:
\begin{equation}
 d_x(\by, \hat{\boldmu}_x)= 2\{\ell_x(\by; \tilde{\btheta}) - \ell_x(\by; \hat{\btheta}_x) \}= 2 \sum_i [ y_i \{ \tilde{\theta}_i - \hat{\theta}_i(x)\} - b( \tilde{\theta}_i) + b( \hat{\theta}_i(x)) ] K_h (x_i -x),
\label{ldev}
\end{equation}
where $\tilde{\btheta}= (\tilde{\theta}_1, \dots, \tilde{\theta}_1)^{\top}$ with  $\tilde{\theta}_i= G(y_i)$, same as those defined around (\ref{paradev}).
Though (\ref{ldev}) is a natural definition for local likelihood, we are not aware of a similar quantity  to (\ref{ldev}) in the literature. Published work focuses on global deviance by taking (\ref{paradev}) with $ G^{-1}\{ \hat{\beta}_0(x_i)\}$ as the estimate, and strictly speaking, the resulting deviance  expression  is not based on maximized likelihood as $\hat{\beta}_1, \dots, \hat{\beta}_p$ are ignored.  
In comparison, the  deviance (\ref{ldev}) makes use of all coefficients $\hat{\beta}_0, \dots, \hat{\beta}_p$ from maximizing local likelihood.
Then (\ref{simon}) is adapted to form a local analysis of deviance expression, and    a global expression may be obtained by integrating local quantities, as given 
 in the following Theorem.

\begin{theorem} 
Suppose that conditions (A1)  and (A2) in the Appendix hold.
Under model (\ref{nonpar}),   
the following results  hold   when using local polynomial approximations of $p$-th order.\\
(a)  For a grid point $x$ in the support of covariate $X$,  a local analysis of deviance expression  is
\begin{equation}
 d_x(\by, \bar{y})= d_x(\by, \hat{\boldmu}_x) + d_x( \hat{\boldmu}_x, \bar{y}),
\label{laod}
\end{equation}
where  
$\bar{y}$ is the sample mean of $\by$, $d_x(\by, \bar{y})$ is  (\ref{ldev}) with $\hat{\boldmu}_x$ and $\hat{\btheta}_x$ replaced by $\bar{y}$ and $G(\bar{y})$ respectively, and  
\begin{equation}
 d_x( \hat{\boldmu}_x, \bar{y}) \equiv 2\E_{\hat{\boldmu}_x} \left[ \ell_x(\by; \hat{\btheta}_x) - \ell_x \{\by; G^{-1}(\bar{y})\} \right] = 2 \left [ \ell_x(\by; \hat{\btheta}_x) - \ell_x \{\by; G^{-1}(\bar{y})\} \right ].
\label{eq3}
  \end{equation}
(b) A global analysis of deviance expression is obtained by integrating the local quantities in (\ref{laod})  over the support of covariate $X$:
\begin{equation}
\int  d_x(\by, \bar{y}) dx= \int  d_x(\by, \hat{\boldmu}_x) dx + \int d_x( \hat{\boldmu}_x, \bar{y}) dx,
\label{gaod}
\end{equation}
where 
 $\int  d_x(\by, \bar{y}) dx=KL(\by, \bar{y})= D(\by, \bar{y})$ under a boundary condition in (A1) that the weights $\int K_h(x_i-x) dx=1$, $i=1, \dots, n$.
\end{theorem}

Theorem 1 provides elegant local and global analysis of deviance expressions that mimic the classic case (\ref{simon}) (McCullagh and Nelder 1989; Simon 1973) and shows that the 
Pythagorean property of the KL divergence holds under model (\ref{nonpar}) with local polynomial fitting. 
It is straightforward to show  (\ref{laod})   based on    (\ref{eq1})  and  (\ref{eq3}) and hence the proof is omitted. Alternatively  the proof in Simon (1973) for (\ref{simon}) can be adapted with kernel weights to show (\ref{laod}). The local expression (\ref{laod}) has an interpretation that  the null deviance at point $x$,  $d_x(\by, \bar{y})$, can be decomposed into two parts, the residual deviance after fitting a locally weighted polynomial at $x$, $d_x(\by, \hat{\boldmu}_x)$, and the model deviance at $x$, $d_x( \hat{\boldmu}_x, \bar{y})$.  Equality (\ref{laod})   holds in finite-sample cases, similar to (\ref{simon}).  The global analysis of deviance    (\ref{gaod})   extends the above interpretation to a fitted curve by local polynomials:
 the residual deviance $ \int  d_x(\by, \hat{\boldmu}_x) dx$ is a measure of the lack of fit of fitting  (\ref{nonpar}), whereas the null deviance  $\int  d_x(\by, \bar{y}) dx$ is such a measure for a reduced model that only includes the intercept.
The  quantities in (\ref{gaod}) are   weighted integrals (see (\ref{lx})), which may be  approximated by the Riemann sum in practice,  and
 an analysis of deviance table based on (\ref{gaod})  is formed, similar to the parametric framework. 
 For a special case of  Normal distribution with an identity link, (\ref{laod}) and (\ref{gaod})  become the local and global analysis of variance decompositions respectively in Huang and Chen (2008).
For the boundary condition in Theorem 1(b),  if $K(\cdot)$ has a support $[-1,1]$ and $\{x_i, i=1, \dots, n\}$ has a range of $[a,b]$, then  a boundary-corrected kernel $[\int K_h(x_i-x) dx]^{-1}K_h(x_i-x) $  may be used for $x_i$ in $[a,a+h)$ and $(b-h, b]$ 
 to ensure that the integrated kernel weights equals to 1.

 As a by-product, the above derivations give rise to new ``global" estimators for $\theta_i$'s and $\mu_i$'s:
\begin{eqnarray}
  \theta_i^*= \int \hat{\theta}_i (x)K_h (x_i -x)  dx \hspace{1cm} \mbox{and} \hspace{1cm}
\mu_i^*= G^{-1}( \theta_i^*).
\label{thetastar}
\end{eqnarray}
They are different from local  estimates at $x_i$:  $\hat{\beta}_0(x_i)$ and $G^{-1} \{ \hat{\beta}_0(x_i)\}$.  The asymptotic properties of $\theta_i^*$ and $\mu_i^*$ for ``interior" points $x_i$ with  $p=1$ and $3$ are discussed in the Proposition  below. The reason  $p=1$ and 3 only is due to their simpler asymptotic bias expressions of $\hat{\beta}_0(x)$ than those of $p=0$ and 2; see Theorems 1a and 1b in Fan et al. (1995). 
The  ``interior"  region is defined as follows. For a kernel function with support  $[-1,1]$, if the convex support of $x_i$s is $[a,b]$, then define the interior region as $[a+2h, b-2h]$.  This definition is narrower than
 the conventional $[a+h, b-h]$, since for $x_i$ in $[a+h, a+2h)\bigcup (b-2h, b-h]$, the corresponding
 $\theta_i^*$ and $\mu_i^*$ in (\ref{thetastar}) involve $\hat{\beta}_j(x)$ with $x$ in $[a, a+h)\bigcup (b-h, b]$, $j=0, \dots, p$.

\begin{proposition}
 Suppose that conditions (A1)-(A5) in the Appendix hold. Assume that  $h \rightarrow 0$ and $nh^3 \rightarrow \infty$ as 
$n \rightarrow \infty$. Then for interior points $x_i$ with  $p=1$ and $3$, \\
  (a) the order of the asymptotic bias of $\theta_i^*$ is smaller than the conventional order $h^{(p+1)}$; i.e.,  the $h^{(p+1)}$ term of the bias of $\theta_i^*$ is zero;\\
(b)  the asymptotic variance of $\theta_i^*$ is of order $n^{-1} h^{-1}$;\\
(c)  similarly,    the order of the  bias of $\mu_i^*$ is smaller than the conventional order $h^{(p+1)}$ and the asymptotic variance of $\mu_i^*$ is of order $n^{-1} h^{-1}$. 
\end{proposition}
The proof for the Proposition  is given in the Appendix. There has been some research aimed at finding new ways of reducing bias of basic kernel smoothers, e.g.,  Kosmidis and Firth, (2009).  
In the Gaussian case with an identity link, Huang and Chan (2014) show that the bias of $\theta_i^*$ for interior points  is of order $h^{2(p+1)}$ for $p=0,1,2,3$, which is consistent with intuition that the higher the $p$, the smaller the order of the bias.  The derivation of explicit bias expressions of $\theta_i^*$ in exponential family is technically challenging, since  the  second-order expansions of the bias of $\hat{\beta}_j(x)$ for  (\ref{glm}) with (\ref{nonpar}) have not been addressed in the literature.
We thus focus on analysis of deviance, while the issue of  bias reduction  may be studied in a future paper.

Theorem 1 involves integrating local likelihood quantities to form a global analysis of deviance expression and 
hence it is of interest to explore how integrated local  likelihood $\int  \ell_x(\by, \hat{\btheta}_x) dx$ behaves as a global likelihood function.
The following theorem shows that  integrated local  likelihood  is  asymptotically a global likelihood $\ell(\by; \btheta^*)$ with estimate  $\btheta^* = (\theta_1^*, \dots, \theta_n^*)^{\top}$  and that the integrated deviance quantities $\int d_x(\by, \hat{\boldmu}_x) dx$ and
$\int d_x( \hat{\boldmu}_x, \bar{y}) dx$ are asymptotically KL-divergence measures with estimate $\boldmu^* = (\mu_1^*, \dots, \mu_n^*)^{\top}$.

\begin{theorem}
Under model (\ref{nonpar}), assume that conditions (A1)-(A5) in the Appendix  hold, and $h \rightarrow 0$, $nh^3 \rightarrow \infty$ as 
$n \rightarrow \infty$.  For  $p=1$ and $3$,\\
(a) the integrated likelihood function is asymptotically
 \begin{equation}
\int  \ell_x(\by; \hat{\btheta}_x) dx= \ell(\by; \btheta^*) + \CO(h^{(p+1)}),
\label{eq2}
\end{equation}
where   the elements of  $\btheta^*$ are defined in (\ref{thetastar});\\
(b) the integrated deviance quantities are asymptotically
 \begin{equation}
 \int  d_x(\by, \hat{\boldmu}_x) dx= KL (\by, \boldmu^*) +\CO(h^{(p+1)}) \hspace{0.5cm}  \mbox{and} \hspace{0.5cm}
\int d_x( \hat{\boldmu}_x, \bar{y}) dx = KL(\boldmu^*, \bar{y}) +\CO(h^{(p+1)}),
\label{KL2}
\end{equation}
where the elements of $\boldmu^*$ are defined in (\ref{thetastar});\\
 (c) from    (\ref{gaod}) and  (\ref{KL2}),
\begin{equation}
  KL(\by, \bar{y})= KL (\by, \boldmu^*)+KL(\boldmu^*, \bar{y})+\CO(h^{(p+1)}),
\label{KLlp}
\end{equation}
which shows that the classic analysis of deviance holds asymptotically with $\boldmu^*$. 
\end{theorem}
The proof of Theorem 2  is given in the Appendix and it utilizes some results stated  in the Proposition  for $p=1$ and 3.  Hence  Theorem 2 is limited to $p=1$ and 3 only for the same reason described before the Proposition. The integrated local likelihood $\int  \ell_x(\by, \hat{\btheta}_x) dx$ in  (\ref{eq2}) is a weighted integral of  local likelihood with fitted local polynomials. 
In the literature,  the idea of integrated likelihood was mentioned  in  Lehmann (2006), and Severini (2007) discussed integrated likelihood functions to eliminate nuisance parameters in parametric settings.  To our knowledge, (\ref{eq2})
and (\ref{KLlp})  have never been raised in the nonparametric regression literature. The convention was to plug in $\hat{\beta}_0(x_i)$ in (\ref{glm}) for $\theta_i$; as $\hat{\beta}_0(x_i)$'s are not maximum likelihood estimates globally, the KL-type additivity (\ref{simon}) would not hold. In contrast,   (\ref{KLlp})
shows that the classic analysis of deviance holds asymptotically by utilizing the local additivity in (\ref{laod}).
Two topics for further investigation are to apply Theorems 1 and 2  to develop bandwidth selection and residual diagnostic procedures. For example,  bandwidth selection may be based on  cross-validating the deviance or minimizing the corrected Akaike information
criterion  (AICc, Hurvich et al. 1998), both with close connection to KL divergence. In Section 5, we explore adapting the 
AICc criterion  with the integrated deviance for  bandwidth selection empirically.

Based on integrated local likelihood,
we next develop an  integrated likelihood ratio test  for examining the significance of a nonparametric fit, parallel to chi-square tests in parametric settings (McCullagh and Nelder 1989). Under model (\ref{nonpar}), the intercept term is embedded in $m(\cdot)$ and hence testing significance of $m(\cdot)$ becomes testing whether $m(\cdot)$ equals to a constant.

\begin{theorem}
Under the conditions of Theorem 2,    for  testing   $H_0: m(x)=a_0$ with $a_0$ a constant  versus $H_a: m(x)$ is not a constant function, when estimating $m(\cdot)$ by $p$-th order local polynomials with $p \geq 0$,  the test statistic
\begin{equation} 
2\left \{\int  \ell_x(\by; \hat{\btheta}_x) dx - \ell(\by; \hat{a}_0) \right\}
\label{liktest}
\end{equation}
is asymptotically distributed according to a $\chi^2$-distribution with  degrees of freedom (df) \text{tr}$(H_p^*) -1$, where $\hat{a}_0$ is the maximum likelihood estimate  under $H_0$ and $H_p^*$ is the smoothing matrix for local $p$-th order polynomial regression defined in Huang and Chen (2008) in the case of the Normal distribution. 
\end{theorem}
More explicitly, 
$H_p^*$  depending on $x_i$s, bandwidth $h$, and the kernel function $K(\cdot)$, is
\begin{equation}
H_p^*=\int W  X_p (X_p^{\top} W X_p)^{-1}X_p^{\top} W dx, 
\label{Hstar}
\end{equation}
where $W$ is an $n$-dimensional  diagonal matrix with $K_h(x_i-x)$ as its diagonal elements, and $X_p$ is the $n \times (p+1)$ design matrix with the $(j+1)$-th column $( (x_1 -x)^{j}, \dots, (x_n-x)^{j})^{\top}$, $j=0, \dots, p$. The  dependence of $W$ and $X_p$ on $x$ is  suppressed and the integration in (\ref{Hstar}) 
is performed element by element in the resulting matrix product. 
In Theorem 3, the $\chi^2$-distribution is allowed to have a non-integer degree of freedom, since the $\chi^2$-distribution is a special case of the gamma
distribution. The asymptotic order of $tr(H_1^*)$ in the case of local linear regression $p=1$ is of order $h^{-1}$ (Huang and Chen 2008,   p. 2093).   
We name the  $\chi^2$-test in Theorem 3 as an  integrated likelihood ratio test since the test statistic can be expressed as integrated likelihood ratio:
\[ \int  \ell_x(\by; \hat{\btheta}_x) dx - \ell(\by; \hat{a}_0)=
\int \sum_i [\ell\{y_i; \hat{\theta}_i(x)\} -  \ell(y_i; \hat{a}_0)] K_h(x_i-x) dx\]
under the boundary condition in (A1).
In other words, under model (\ref{nonpar}), the test statistic (\ref{liktest})  integrates the differences in local deviances between a nonparametric fit (\ref{ldev})   and an intercept-only reduced model and  it is distributed asymptotically as a chi-squared distribution with the difference in degrees of freedom of the two models. This interpretation makes 
(\ref{liktest}) more compelling than the generalized likelihood ratio test in Li and Liang (2008), since their work does not have a connection to deviance.

\section{Analysis of Deviance for Partially Linear Models}

We extend the results in Section 3 to  generalized partially linear models (\ref{model}).
Denote $\breve{\mu}_x(x_i)= G^{-1} \{\breve{\theta}_i(x)\}$ where
$\breve{\theta}_i(x)= \bz_i^{\top} \breve{\alpha} + \breve{\beta}_0 + \dots + \breve{\beta}_p (x_i -x)^p$ with estimates $\breve{\alpha}$ and $\breve{\beta}_j$'s under (\ref{model}). To avoid confusion with the notation in Section 3, from now on $\breve{\mu}_x(x_i)$, $\breve{\theta}_i(x)$, $\breve{\btheta}_x$, $\breve{\boldmu}_x$, $\breve{\alpha}$, $\boldmu^{**}$, and $\btheta^{**}$  denote the estimates under 
 (\ref{model}).
Since $\breve{\beta}_j$'s maximize the local likelihood,  the  equations in 
(\ref{eq1}) continue to  hold with  $\breve{\btheta}_x$ and $\breve{\boldmu}_x$ under (\ref{model}). The interpretation that  $(\by - \breve{\boldmu}_x)$ is orthogonal to 
$\breve{\btheta}_x$ in the locally weighted inner product space with weights 
$K_h (x_i -x)$ continues to hold under (\ref{model}). An additional equation from  estimating $\alpha$ by maximum likelihood is
\begin{equation}
  \sum_i y_i z_{ik}     =  \sum_i z_{ik} \int \breve{\mu}_x(x_i) K_h(x_i-x) dx, \hspace{1cm} k = 1, \dots, K,\label{eq6}  
\end{equation}
where $z_{ik}$ denotes  the value of the $k$-th covariate for the $i$-th observation.  
From  (\ref{eq6}), we observe that the column vector with entries  $(y_i - \int \breve{\mu}_x(x_i) K_h(x_i-x) dx)$, $i=1, \dots, n$, is  orthogonal to the column space spanned by $\bz$.
Moreover it can be shown that $\int \breve{\mu}_x(x_i) K_h(x_i-x) dx=\mu_i^{**}+ \CO(h^{p+1}+n^{-1/2})$ and 
 hence $(\by - \boldmu^{**})$ is asymptotically orthogonal to the column space spanned by $\bz$.

Theorems 1 and 2  are extended  to generalized partially linear models (\ref{model})  in the following as Theorem 4(a) and 4(b) respectively 
 when $\alpha$ is estimated by maximum likelihood. We develop local and global analysis of deviance expressions for (\ref{model}) in Theorem 4(a)  and Theorem 4(b) shows that
 the integrated likelihood quantities are asymptotically global likelihood quantities with $\btheta^{**}$  and  $\boldmu^{**}$.  

\begin{theorem}  For model   (\ref{model}), assume that Conditions (A) in the Appendix  hold,  and $h \rightarrow 0$, $nh^3 \rightarrow \infty$ as 
$n \rightarrow \infty$.  \\
(a)  The local and global analysis of deviance (\ref{laod}) and (\ref{gaod}) respectively hold  with $\breve{\mu}_x(x_i)$ and $\breve{\boldmu}_x$ when $\alpha$ is estimated by maximum likelihood.\\
(b) Assume that $\alpha$ is estimated with a root-$n$ rate. For $p=1$ or 3, the  expression in  (\ref{eq2}) holds with  $\breve{\btheta}_x$ and $\btheta^{**}$ except the $\CO(h^{p+1})$ term replaced by
  $\CO(h^{p+1}+n^{-1/2})$.  Similarly,  (\ref{KL2})  and (\ref{KLlp}) hold  with $\breve{\mu}_x$ and $\boldmu^{**}$ and the $\CO(h^{p+1})$ terms replaced by
  $\CO(h^{p+1}+n^{-1/2})$. \\
(c) When the same kernel function and bandwidth are used in  (\ref{nonpar})
and (\ref{model}),
 the nonparametric model (\ref{nonpar}) is nested in 
(\ref{model}). Then  the difference in  local residual deviance from fitting (\ref{model}) to  (\ref{nonpar})  can be expressed  as
\begin{equation}
d_x (\by, \hat{\boldmu}_{x} ) - d_x(\by, \breve{\boldmu}_x) = d_x(\breve{\boldmu}_{x}, 
 \hat{\boldmu}_{x}).
\label{aod2}
\end{equation}
\end{theorem}
The proofs of Theorem 4(a) and  4(b) are analogous to Theorems 1 and  2 respectively and are thus omitted. 
We briefly outline the proof for Theorem 4(c).   Based on   
(\ref{eq1}) under (\ref{model}), we have $\sum_i y_i (x_i -x)^j  K_h (x_i -x)  =  \sum_i \breve{\mu}_x(x_i)(x_i -x)^jK_h (x_i -x),  j = 0, \dots, p$. Then multiplying  the $j$-th equation  by 
$\hat{\beta}_j$  and summing them up, 
 $\sum_i \{ y_i - \breve{\mu}_x(X_i)) \} \hat{\theta}_{i}(x) K_h (X_i -x)=0$ is obtained and  (\ref{aod2}) is proved.

In a special case of  the Gaussian distribution with an identity link, Theorem 4(a)  becomes the local and global analysis of variance  for partially linear models, which was discussed in Huang and Davidson (2010, section 3.1).
Theorem 4(c) implies that  the local residual deviance for fitting (\ref{model}) is the local residual  deviance for fitting   (\ref{nonpar}) minus a term due to   the parametric component. That is,  the difference of local  residual deviances between  (\ref{model}) and (\ref{nonpar}) is a KL-divergence measure $d_x(\breve{\boldmu}_{x}, 
 \hat{\boldmu}_{x})$, and the local KL-divergence  is additive between nested models  (\ref{model}) and (\ref{nonpar}). A similar interpretation holds at a global scale after integrating the local contributions of (\ref{aod2}):
\[
\int d_x(\by, \breve{\boldmu}_x)  dx = \int d_x (\by, \hat{\boldmu}_{x} ) dx - \int d_x(\breve{\boldmu}_{x},  \hat{\boldmu}_{x}) dx.
\]

Analogous to Theorem 3,  testing whether $m(\cdot)$ is significant under model   (\ref{model}) becomes 
testing whether $m(\cdot)$ is significantly different from a constant and an  integrated  likelihood ratio test is proposed in the following theorem.

\begin{theorem} 
 For model   (\ref{model}), assume that Conditions (A) in the Appendix  hold and  the data matrix $Z$ for covariates $\bz$ is orthogonal to $\bx=(x_1, \dots, x_n)^{\top}$ and the intercept column.
For  testing   $H_0: m(x)=a_0$ with $a_0$ a constant  versus $H_a: m(x)$ is not a constant function, when  estimating $m(\cdot)$ by $p$-th order local polynomials with $p \geq 0$,  the test statistic
\begin{equation} 
2\left \{\int  \ell_x(\by; \breve{\btheta}_x) dx - \ell(\by; \hat{\balpha}_0) \right\}
\label{liktest2}
\end{equation}
is asymptotically distributed according to a $\chi^2$-distribution with  df $\text{tr}(H_p^*) - 1$, where  $\hat{\balpha}_0$
is the maximum likelihood estimate  under parametric $H_0$.
\end{theorem}

The assumption  that $Z$ is orthogonal to $\bx$  in Theorem 5 is required for mathematical convenience, in the sense that the corresponding off-diagonal elements of the local Fisher information is 0, for ease of deriving  the asymptotic $\chi^2$-distribution of the test statistic (\ref{liktest2}). It will be seen in the simulations that the performance of the proposed tests remain reasonable when this assumption is violated.
The integrated likelihood ratio tests  in Theorems 3 and 5 depend on the bandwidth $h$, like the other nonparametric tests. Analogous to Theorem 3, the test statistic (\ref{liktest2}) has an interpretation of  integrating the differences in local deviances between a fitted generalized partially linear model   and an intercept-only reduced model. 
We remark that the proposed tests are different from those in Hastie and Tibshirani (1990) and Li and Liang (2008). The proposed test statistics utilize  integrated  likelihood that combines all maximized local likelihoods from fitting local polynomials. Some existing methods use only  $G^{-1}(\hat{\beta}_0(x_i))$'s, and strictly speaking, the resulting expression  is not based on maximizing likelihood as $\hat{\beta}_1, \dots, \hat{\beta}_p$ are ignored; this fact was mentioned before Theorem 1 as well.  

Some work in the literature, e.g. H\"ardle et al. (1998),  has considered testing whether $m$ in (\ref{model}) is significantly different from a linear trend,  $G(\mu)= \bz^{\top} \alpha + a_0+a_1 x$.  The extension of Theorem 5 to testing a linear trend is non-trivial and will be pursued in future work,
since the variance function of $y$ is allowed to be a function of the mean of $y$ in exponential family (\ref{glm}).
In a special case of the Gaussian distribution with a constant variance,   analysis-of-variance $F$-type tests for checking linear  trends are derived in Huang and Davidson (2010).

\section{Simulation Results}

We  examine the empirical type-I errors and power for the proposed tests  in Theorems 3 and 5. Local linear smoothing $p=1$ with  the Epanechnikov kernel  is used throughout this section.  We first describe the algorithm for calculating  the test statistic (\ref{liktest2}) in Theorem 5  for testing $H_0: m(x)=a_0$ under model   (\ref{model}),  while that for  (\ref{liktest}) under model (\ref{nonpar})   is similar. 
The algorithm   adapted from  Carroll et al. (1997) 
 is given as follows:\\
\hspace{1cm} Step 0 (initialization). Fit a parametric generalized linear model to obtain initial values $\breve{\alpha}^{(0)}$.\\
\hspace{1cm} Step 1. For a set of grid points  on the data range of $x_i$'s,  given a value of $h$,   maximize the local likelihood with  $\breve{\alpha}^{(r)}$ to obtain $\breve{\beta}_0^{(r)}, \dots, \breve{\beta}_p^{(r)}$ for each grid point. Then with $\breve{\theta}_i^{(r)}(x)$'s,  calculate a locally weighted average as in  (\ref{thetastar}) to obtain $\breve{\theta}_i^{**^{(r)}}$, $i=1, \dots, n$. \\
\hspace{1cm} Step 2.   Maximize the  global likelihood with $\btheta^{**^{(r)}}=(\breve{\theta}_1^{**^{(r)}}, \dots, \breve{\theta}_n^{**^{(r)}})^{\top}$ to update $\breve{\alpha}^{(r+1)}$.\\
\hspace{1cm} Step 3. Continue Steps 1 and 2 until convergence. The test statistic (\ref{liktest2}) is calculated by integrating the final  local likelihoods and taking its difference to the global likelihood under $H_0$.

The simulation study focuses on logistic regression in Examples 1-4 as we wish to evaluate the proposed methods in order to analyze the German Bundesbank data in Section 6, while Example 5 is on Poisson regression. 
The integrated likelihood in the test statistics (\ref{liktest}) and (\ref{liktest2}) are approximated discretely 
by taking 201 equally-spaced points on $[0, 1]$ in Examples 1 and 2, and 301 equally-spaced points on $[-0.5, 1]$ in Examples 3-5. For $x_i$s that fall in conventional boundary regions, analogous approximations are used for calculating 
 $\int K_h(x_i-x) dx$  for    boundary correction      in condition (A1).
In addition to implementing the proposed tests with a fixed $h$, 
we also try selecting the bandwidth by AICc (Hurvich et al., 1998)
and  by the idea of Horowitz and Spokoiny (2001) (HS). 
The AICc criterion is adapted with   the integrated deviance: 
\[ AICc(h)= \log(D^*/ n) + 2(tr(H_1^*) +1) /(n - tr(H_1^*) -2),\]
where $D^*$ denotes the integrated deviance $\int  d_x(\by, \hat{\boldmu}_x) dx$ in (\ref{KL2}) under model (\ref{nonpar}) or $\int  d_x(\by, \check{\boldmu}_x) dx$ under model (\ref{model}).
The HS idea  is to select the bandwidth that maximizes the test statistic. 
Critical values  for the proposed tests are taken  from
the  $\chi^2$-distribution with 5\% significance level and   5000 simulated data sets are generated.  The {\tt gam} function in the {\tt mgcv} R-package (Wood,  2013) provides a  chi-square test of zero effect of a smooth term and we include it for comparison, with  default 10 spline basis functions and the penalty estimated by REML.

 {\it Example 1:} $ \text{logit}(p)= -1 + a \cos(2\pi x)$, $a=0, 0.5, 0.75, 1$.
We first check the $\chi^2$-approximation under $H_0$ when $a=0$ for both a fixed design, $x$ equally-spaced on $[0,1]$,  and a random design,  $x\sim U(0,1)$, with sample sizes
$n=$50, 100, and 200. The values of bandwidth $h=$0.1, 0.12, 0.15, 0.17, 0.2,  0.25, and 0.3 are chosen so that they are roughly equally-spaced on a logarithm scale and they correspond to smoothing with about 20\%-60\% data. For both AICc and HS, the bandwidth among the 7 values that satisfies the criterion is selected. The results with  $h=$0.1,  0.15, 0.2,  and 0.25,  are chosen to present in Table 1, from under-smoothing slightly to over-smoothing slightly. The results with varying $h$ by   AICc and HS are also given in Table 1.
It appears that  when $n=$50, the  $\chi^2$-approximation for  (\ref{liktest})  under $H_0$ is not good  as the empirical type-I errors  are all above 0.05 for either a  fixed or random design.  For this reason, we do not consider the case of $n=$50 further.   As suggested by two reviewers, a bootstrap alternative for calculating the sample critical values for $n=50$ may be considered for future research.

\begin{table}
\begin{center}
\caption{Percent of rejection under $H_0$ in Example 1 with $a=0$}
\caption*{\footnotesize  The empirical type-I errors of the test statistic (\ref{liktest})   are  close to 0.05 for $n=100$ and $200$ with a fixed $h\geq0.15$. The performance of 
$h_{AICc}$ when $n=200$ is closer to $H_0$ than that of $n=100$, and 
 $h_{HS}$  has  larger type-I errors as it attempts to optimize the power. The {\tt gam} function performs consistently around 0.05 under $H_0$ regardless of the  sample sizes.}
\begin{tabular}{lccccccc}
 \hline \hline
& $h=0.1$ & $h=0.15$ & $h=0.2$ & $h=0.25$ & $h_{AICc}$ & $h_{HS}$  &  {\tt gam} \\   
$n=50$ fixed design & 12.80 & 8.74 & 6.70 & 6.02 & 10.02 & 15.86 & 3.22 \\ 
$n=50$ random design & 14.10 & 9.34 & 7.20 & 6.24 &10.38 & 17.26 &  3.06 \\  \hline
$n=100$ fixed design &  ~7.38 & 5.30 & 4.80 & 4.52 & ~7.18 & ~9.96 &  3.99 \\
$n=100$ random design &    ~7.94 & 5.84 & 4.70 &  4.36 &  ~7.28 & 10.48 & 4.08  \\  \hline
$n=200$ fixed design &  ~4.96 & 4.20 & 4.18 & 3.96 & ~5.56& ~7.36 &  4.62 \\
$n=200$ random design & 5.22 & 4.38 & 4.42 & 4.20 & ~6.00 & ~7.76 & 5.04 \\ \hline \hline
\end{tabular}
\end{center}
\end{table}

When $a=0$ and $n=100$, the empirical type-I errors are mostly reasonable except for a small $h=$0.1,  $h_{AICc}$, and $h_{HS}$, with  rates ranging  about $7-10\%$. In this case with a random design,  $h_{AICc}$  tends to select the largest bandwidth 0.3,   92.12\% of 5000 simulations, since the true model under $H_0$ is a constant, and 
when AICc happens to select a small bandwidth such as 0.1, it often leads to rejecting $H_0$.  
For $h_{HS}$, it  behaves differently since it attempts to optimize the power; when $a=0$ and $n=100$, the empirical proportions of $h_{HS}$ on the 7 values of $h=0.1, \dots, 0.3$ are  44.58\%, 5.8\%, 4.96\%,  4.08\%,  5.04\%,  4.84\%, and  30.70\%.
Therefore the inflated type-I errors of   $h_{HS}$ are somewhat  expected. In Horowitz and Spokoiny (2001), the critical values was based on resampling from the finite-sample null distribution, while we use the asymptotic $\chi^2$-distribution.  
When $n=200$, the empirical type-I errors for the proposed tests are around $0.05$ with a fixed bandwidth, and slightly above 0.05 for $h_{AICc}$ and  $h_{HS}$. The performance of 
$h_{AICc}$ when $n=200$ is closer to $H_0$ than that of $n=100$.
The {\tt gam} function performs consistently around 0.05 under $H_0$ regardless of the  sample sizes.
When $n=$100 with a fixed design, the  df ($tr(H_1^*)-1)$  are 10.29, 6.84, 5.11, and 4.07  respectively for  $h=$0.1, 0.15, 0.2, and 0.25 respectively and  the average estimated degrees of freedom  (edf)   for {\tt gam} is 1.33 with a range $[1.00, 6.85]$. Quantile-quantile plots (qqplots) of 5000 test statistics (\ref{liktest}) for $n=$100 with a fixed design   against the $\chi^2$-quantiles with the corresponding df are shown in Figure 1 for  $h=$0.1, 0.15, 0.2, and 0.25,  indicating  satisfactory approximations of the $\chi^2$-distribution. The qqplots of $n=$200 with a fixed bandwidth  (not shown) are similar to those of $n=$100. 

\begin{figure}
\includegraphics[width=6.5in,height=8in]{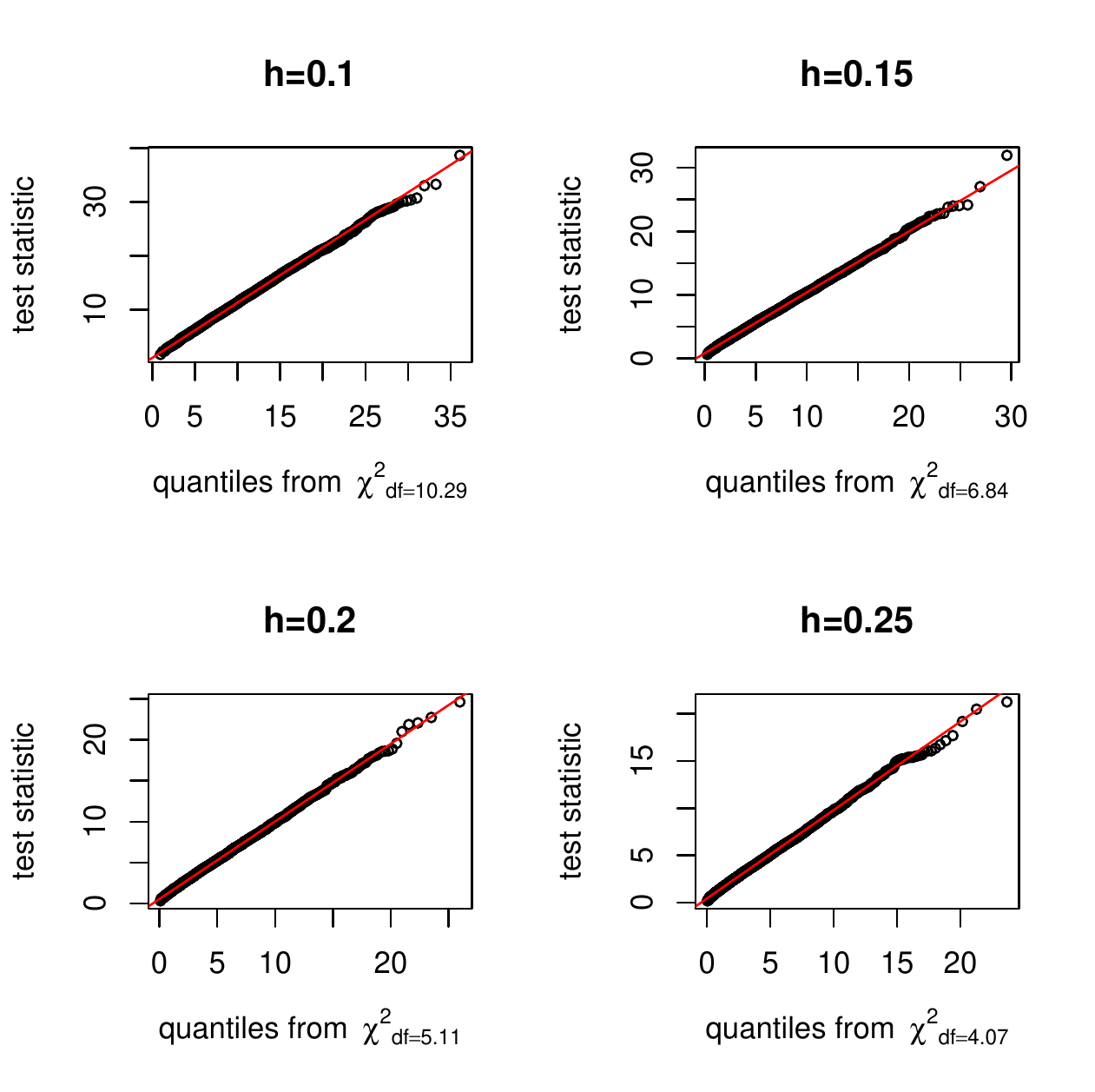} 
\caption{Quantile-quantile plots of  5000 integrated likelihood ratio test  statistics (\ref{liktest}) under $H_0$ in Example 1 for $n=100$  with $h=$0.1, 0.15, 0.2, and 0.25, against quantiles from $\chi^2$-distribution with df 10.29, 6.84, 5.11, and 4.07 respectively.}
\end{figure}

When $a=$0.5, 0.75, and 1, with a random design, we examine the performance under  alternatives. 
The percent of rejection  is given in Table 2.  
For the proposed tests, we observe that the rejection rate increases as the value of the bandwidth increases  when $a=$0.75 and 1, and $h_{AICc}$ and $h_{HS}$ are more powerful than those with a fixed $h$. Under alternatives, the proposed tests are more powerful than {\tt gam} except the case with $n=200$ and $h=$0.1. 
When $n=$100,  the average dfs  for  (\ref{liktest})  of $a=$ 0.5, 0.75, and 1 are similar to those of $a=$0,  since $H_1^*$ in (\ref{Hstar}) does not involve the response $y$. When $n=100$, the average edf  of  {\tt gam} increases as $a$ increases,   1.80, 2.30, and  2.64 for $a=$0.5, 0.75, and 1 respectively.  The behaviour for df and edf of $n=$200 is similar to that of $n=100$.

\begin{table}[htbp]
\begin{center}
\caption{Percent of rejection under $H_1$ in Example 1}
\caption*{\footnotesize The rejection rate for  the test statistic (\ref{liktest})  increases as the value of the bandwidth increases  when $a=$0.75 and 1, and $h_{AICc}$ and $h_{HS}$ are more powerful than those with a fixed $h$. The  test statistic (\ref{liktest})  is more powerful than {\tt gam} except the case with $n=200$ and $h=$0.1.}
\begin{tabular} {lccccccc}
 \hline \hline
 & $h=0.1$ & $h=0.15$ & $h=0.2$ & $h=0.25$ & $h_{AICc}$ & $h_{HS}$  & {\tt gam} \\  
$n=100$, $a=0.5$ & 20.12 &   18.78 &  19.32 &  20.38  & 24.26 & 28.58 & 13.96\\ 
$n=100$, $a=0.75$ & 37.22 &  38.38 &  40.96 &  43.42 & 47.14 & 51.26 & 31.96 \\ 
$n=100$, $a=1$ &  60.56 & 63.72 & 67.96 & 70.60 & 73.52 & 75.78 &  58.86 \\   \hline
$n=200$, $a=0.5$ & 29.18 & 33.04 & 36.30 & 38.96 & 42.28 & 45.08 & 30.82\\
$n=200$, $a=0.75$ & 61.18 & 67.40 & 71.68 & 75.00  & 77.66 & 78.84 &  66.58\\ 
$n=200$, $a=1$ & 88.58 & 92.38 & 94.28 & 95.56  & 96.34 & 96.64 & 92.36 \\   \hline  \hline
\end{tabular}
\end{center}
 \end{table}

{\it Example 2:} $ \text{logit}(p)= -2 + f_k(x)$,    $k=0, 1, 2,$
where   $x\sim U(0,1)$ and  functions $f_0(x)=8x(1-x)$, $f_1(x) =\exp(2x)$, and $f_2(x) =2 \times10^5x^{11}(1-x)^6 + 10^4 x^3 (1-x)^{10}$, are taken from Wood (2013). Wood (2013) considered an additive model with  $\text{logit}(p)= -5+ f_0(x_0) + f_1(x_1) + f_2(x_2)$, while we use those functions in the univariate case separately. 
The results are shown in Table 3, with the proportions of rejection nearly 100\% for both tests in  cases of  $f_1$ and $f_2$. For $f_0$, a quadratic trend, our test is more powerful than 
{\tt gam}  except the case  with $n=200$ and $h=$0.1.
  
\begin{table}[htbp]
\begin{center}
\caption{Percent of rejection   for Example 2}
\caption*{\footnotesize The  rejection rates are nearly 100\% for (\ref{liktest})  and {\tt gam} in  cases of  $f_1$ and $f_2$. For $f_0$, a quadratic trend,  (\ref{liktest}) is more powerful than 
{\tt gam}  except the case  with $n=200$ and $h=$0.1.}
\begin{tabular} {lccccccc} 
 \hline  \hline
& $h=0.1$ & $h=0.15$ & $h=0.2$ & $h=0.25$  & $h_{AICc}$ & $h_{HS}$ & {\tt gam} \\  
$n=100$, $f_0$ &   ~47.94 &   ~51.26  &   ~54.72 &   ~57.46  & ~60.20 & ~64.38 &  ~38.70\\ 
$n=100$, $f_1$ &  ~99.86 &   100 &  100 &  100  & 100 & 100 &  ~98.90 \\ 
$n=100$, $f_2$ &  100 & 100 & 100 & 100 & 100 & 100 & ~99.74 \\  \hline
$n=200$, $f_0$ &  ~76.10 &  ~82.00  &  ~85.38 &  ~87.84 & ~89.34 & ~90.18 &   ~80.82\\ 
$n=200$, $f_1$ & 100 & 100 & 100 & 100  & 100 & 100 &  100\\ 
$n=200$, $f_2$ & 100 & 100 & 100 & 100 & 100 & 100 & 100 \\  \hline  \hline
\end{tabular}
\end{center}
 \end{table}

{\it Example 3:}
 $\text{logit}(p)= b_1 z_1 + b_2 z_2+ a \exp(-16 x^2),$ $a=0,1,2,3$, 
where  $z_1$  is  first generated as binary taking values $-$1 and 1 with equal probabilities,  $z_2$ and $x$  are
first generated from a bivariate normal distribution with mean 0, variances 0.5 and 1 respectively, and correlation 0.3.
Then $x$ is transformed to have a uniform distribution on $(-0.5,1).$ To satisfy the conditions in Theorem 5,  $z_1$ and $z_2$ are then made orthogonal to $\bx$ and the intercept vector. 
After the orthogonized  $z_1$ and $z_2$  are obtained, $b_1=0.1$, $b_2=-0.1$.  
To understand how restrictive    the orthogonality assumption   in Theorem 5 is, we  also examine  the performance of (\ref{liktest2}) with the original non-orthogonalized values of  $z_1$ and $z_2$ and same values of    $b_1$ and $b_2$.  

\begin{table}[htbp]
\begin{center}
\caption{Percent of rejection for Example 3}
\caption*{\footnotesize  When $a=0$ with a fixed $h>=0.2$, the empirical type-I errors of (\ref{liktest2}) are close to 0.05.  
When  $a=1$ and 2, (\ref{liktest2})  is more powerful than  the {\tt gam} test, while for $a=3$, the performance of the two tests are close. Under alternatives,
$h_{HS}$ is the most powerful, while $h_{AICc}$ also performs well. The rejection rates for (\ref{liktest2}) are quite close whether $Z$ and $\bx$ are orthogonal or not (the non-orthogonalized version in brackets). }
\begin{tabular} {lcccccccc}
 \hline \hline
 &   $h=0.2$ & $h=0.25$ & $h=0.3$ & $h=0.4$ &  $h_{AICc}$ & $h_{HS}$ & {\tt gam} \\
$n=100$, $a=0$   & ~6.40  & ~5.78& ~5.34  & ~4.96&
~8.10  & 11.04  &  ~4.26 \\ 
  &  [~6.22] &  [~5.64]  &  [~4.98]  & [~4.74] & [~7.70] & [10.46] &  [~4.46] \\ \hline
$n=100$, $a=1$ & 20.50 &  20.54  &  20.64 & 20.76
 & 23.24  & 28.40 &  13.36  \\
 & [20.98]  & [20.32]  & [20.36]  &  [20.12] & [22.46] &   [28.24] & [13.04] \\ \hline
$n=100$, $a=2$   & 69.14 & 70.76  &  72.00  & 72.30 &
 73.50  & 77.14 & 57.72   \\
   &  [68.36] & [70.26]  & [71.54]  &   [71.74]  & [72.84] & [76.62] & [57.08] \\ \hline
$n=100$, $a=3$    & 96.00  & 96.92 & 97.38 & 97.42 & 
97.54 & 98.22 & 92.34 \\ 
  & [95.50] & [96.56] & [96.88]  & [97.00]& [97.16] & [97.86] & [91.46]\\ \hline \hline
$n=200$, $a=0$ & ~5.56  & ~5.24& ~5.20  & ~4.78  &
6.80  & 8.88   & ~5.14  \\
 & [~5.62]  &   [~5.08]  & [~4.90]  &  [~4.62]  & [6.60] & [8.62] & [~5.00] \\ \hline
$n=200$, $a=1$  & 34.86  & 36.94   & 38.04  & 38.46 &
39.94  & 44.24 & 29.26 \\ 
  & [34.22] & [36.06] & [37.58]  &  [38.48]  & [39.94] & [43.76]  & [29.32]  \\ \hline
$n=200$, $a=2$   & 95.92 & 96.58   & 96.92   & 97.16 &
97.36 & 97.70  & 93.50  \\
  & [95.68] & [96.56] & [97.06] &  [97.20] &  [97.32] & [97.54] & [93.26] \\ \hline
$n=200$, $a=3$   & 100 & 100   & 100  & 100 & 
100 & 100 &  99.96  \\  
 & [100] & [100] & [100] & [100] & [100]  & [100] & [99.98]  \\  \hline \hline
\end{tabular}
\end{center}
 \end{table}

The values of bandwidth are 0.15,  0.2,  0.25,     0.3,  and 0.4, so that they are roughly equally-spaced 
on a logarithm scale. The percent of rejection is given in Table 4 for $h=$0.2, 0.25, 0.3, 0.4, $h_{AICc}$, and $h_{HS}$, with the non-orthogonalized version in brackets. The case of $h=0.15$ is not presented due to its inflated type-I errors: 
when $a=0$, $n=100$, and  $h=0.15$, the percent of rejection is   8.72 and 8.28 for  the orthogonalized and non-orthogonalized version respectively. From Table 4, we observe that when $a=0$, the empirical type-I errors are reasonable except 
$h_{AICc}$ and $h_{HS}$. 
Together with the observations in Example 1 under $H_0$, we may imply that optimizing the bandwidth by some criterion may lead to inflated type-I errors for our test in the case of logistic regression.  
When  $a=1$ and 2, our  test is more powerful than  the {\tt gam} test, while for $a=3$, the performance of the two tests are close. 
 Under alternatives,
$h_{HS}$ is the most powerful, while $h_{AICc}$ also performs well, supporting $AICc$ as a bandwidth-selection criterion.  The rejection rates for (\ref{liktest2}) are quite close whether $Z$ and $\bx$ are orthogonal or not, suggesting that this assumption may be relaxed in practice.  When $n=200$ and $a=0$, the average df $(tr(H_1^*)-1)$ corresponding to  $h=0.2, \dots,$ 0.4  are 7.69,  6.14,  5.10, and 4.81 respectively, and again they stay about the same between different values of $a$. 
When $n=200$, the average edf for {\tt gam} is  1.34, 2.22, 3.89, and 4.73 for $a=$0, 1,  2, 3 respectively.  The df and edf of $n=100$ are similar to those of $n=200$.

{\it Example 4:} $\text{logit}(p)= b_1 z_1 + b_2 z_2+ a \cos(2\pi x),$   $a=0.5, 1,  1.5$,
where the data generation scheme of $z_1$, $z_2$, and $x$ is identical  to Example 3, and $b_1$ and $b_2$ are the same as Example 3. This example adopts a  nonlinear function of $x$ similar to that of Example 1 with a different range of $x$. The same values of $h$ as Example 3 are used and hence the dfs are analogous to Example 3, omitted for brevity. Table 5 shows that  (\ref{liktest2}) is more powerful than {\tt gam} when $a=0.5$ and 1.0, and our test with  $h_{HS}$ and $h_{AICc}$ continues to perform well in this example. Again, we observe that for the proposed tests, the rejection rates are quite close whether $Z$ and $\bx$ are orthogonal or not.

\begin{table}[htbp]
\begin{center}
\caption{Percent of rejection for Example 4}
\caption*{\footnotesize   The test statistic (\ref{liktest2}) is more powerful than {\tt gam} under alternatives. The rejection rates are quite close whether $Z$ and $\bx$ are orthogonal or not (the non-orthogonalized version in brackets).}
\begin{tabular} {lccccccc} 
 \hline \hline
 & $h=0.2$ & $h=0.25$ & $h=0.3$ &  $h=0.4$ & $h_{AICc}$ & $h_{HS}$ & {\tt gam} \\ 
$n=100$, $a=0.5$  & 20.16  & 19.66  &  19.62 &18.56 &
22.38  & 26.78 &  ~8.60  \\
  & [20.02] & [19.54] &  [19.22]  & [18.60]& [22.40] & [27.24]  & [~8.98] \\ \hline
$n=100$, $a=1$  & 62.76 & 64.36  & 65.46 &  64.68 &
67.08  & 71.54  & 38.86 \\
 & [61.60]  & [63.48] & [64.08] & [63.38]  & [65.88] & [70.58] & [38.34] \\ \hline
$n=100$, $a=1.5$   & 94.58 & 95.18   &95.60  & 95.18 &
95.72 & 96.48  &   85.24\\ 
  &   [94.00]  & [94.98] & [95.22] & [94.72] &  [95.18] & [96.14] &  [84.90] \\ \hline \hline
$n=200$, $a=0.5$  & 30.60 & 32.36  & 33.36   & 32.58 &
34.90& 39.00  &  18.58 \\ 
 &  [30.80] & [32.24] & [32.88] &  [32.54]  &  [34.94]  & [38.88] &  [18.66]  \\ \hline
$n=200$, $a=1$   & 91.60  & 93.10 & 93.92  & 93.68 &
94.18  & 95.04 &    83.60 \\ 
  & [90.94] & [92.42]  & [93.50] & [93.04] & [93.48] & [94.48]  &  [82.84]  \\ \hline
$n=200$, $a=1.5 $ & 100  & 100  & 100 & 100 & 
100  & 100 &  99.86  \\  
 & [99.98] & [100] & [100] & [100]  & [100] & [100] & [99.88] \\ \hline \hline
\end{tabular}
\end{center}
 \end{table}

 {\it Example 5:}  
\begin{equation}
\text{log}(\mu)= b_1 z_1 + b_2 z_2+ a \exp(-16 x^2), a=0, 1, 2,
\label{eg51}
\end{equation}
 and
\begin{equation}
\text{log}(\mu)= b_1 z_1 + b_2 z_2+ a \cos(2\pi x),  a=0.5, 1.5.
\label{eg52}
\end{equation}
This example is for Poisson regression with the canonical log link, while the  functional form for $\theta=\text{log}(\mu)$ and data generation scheme follow those of Examples 3 and 4. From Table 6,  we observe that  when $a=0$  in (\ref{eg51}) with $n=100$ and $h=0.15$, the type-I error is reasonable, in contrast  to the logistic regression case.   The performance of $h=0.4$ is close to that of $h=0.3$ and therefore  not presented  in Table 6. We observe that our test  using a fixed $h$ is more powerful with a larger bandwidth when $a=1$ in
(\ref{eg51}) and $a=0.5$ in (\ref{eg52}),  and the empirical power is comparable between (\ref{liktest2}) and {\tt gam} tests. 

\begin{table}[htbp]
\begin{center}
\caption{Percent of rejection for Example 5}
\caption*{\footnotesize When $a=0$, the type-I errors of (\ref{liktest2}) are reasonable for this Poisson regression example. 
 The  test statistic (\ref{liktest2})   using a fixed $h$ is more powerful with a larger bandwidth when $a=1$ in
(\ref{eg51}) and $a=0.5$ in (\ref{eg52}),  and the empirical power is comparable between (\ref{liktest2}) and {\tt gam} tests. 
}
\begin{tabular} {lccccccc}
 \hline \hline
 &  $h=0.15$ & $h=0.2$ & $h=0.25$ & $h=0.3$ & $h_{AICc}$ & $h_{HS}$ & {\tt gam} \\ 
$n=100$, $a=0$  in (\ref{eg51}) &  ~5.36  &  ~4.80 & ~4.18  & ~3.94&
~7.06& ~7.62 &  ~4.88 \\ 
 & [~5.48] & [~4.64]  & [~4.40] &   [~4.26]  &  [~7.14]  &  [~7.68] & [~4.86]  \\ \hline
$n=100$, $a=1$ in (\ref{eg51}) &  41.26 &  ~44.48  &  47.24 &  48.42   &
51.98  & 54.48  &  37.58 \\ 
 & [41.36]  & [44.52] &  [47.12]  & [~48.90] & [52.60] & [54.98] & [37.24]  \\ \hline
$n=100$, $a=2$ in (\ref{eg51}) & 99.98 & 99.98& 99.98 & 100 &
99.98 & 100 & 99.94 \\
 & [99.98]  &  [99.98]  & [99.98]  & [100]  & [99.98]  & [100]  & [99.96]\\ \hline
$n=200$, $a=0$  in (\ref{eg51}) & ~4.98 & ~4.56 & ~4.48  & ~4.32 &
~6.92  & ~7.38 & ~5.14 \\
 & [~5.02]  & [~4.66]  & [~4.74]  & [~4.60] & [~7.06] & [~7.62]  & [~5.32] \\ \hline
$n=200$, $a=1$  in (\ref{eg51}) & 74.10 & 78.64   & 81.48  & 82.88 & 85.02 & 85.86 &  75.04   \\
 &  [73.62] & [78.46] & [80.84] & [82.04] & [84.12] & [84.92] & [74.14] \\ \hline
$n=200$, $a=2$  in (\ref{eg51}) & 100 & 100 & 100 & 100 & 100 & 100 & 100 \\
 & [100] & [100] & [100] & [100] & [100] & [100] & [100] \\ \hline \hline
$n=100$, $a=0.5$ in (\ref{eg52}) &  26.20   &  28.36  & 30.02   &  30.90  &
34.72 & 36.58  &  17.32   \\  
& [25.72] & [27.56] & [29.40] & [30.18] & [34.06] & [35.66] & [16.96] \\ \hline 
$n=100$, $a=1.5$ in (\ref{eg52}) & 99.92 & 100 & 100 & 100 & 100 & 100 & 99.98 \\
 & [99.94] & [99.94] & [99.98] & [100] & [100] & [100] & [99.94] \\ \hline 
$n=200$, $a=0.5$ in (\ref{eg52}) & 51.86   &  56.74 & 60.22  & 62.24  & 
64.90  & 65.88  &  43.94 \\
 & [51.18] & [56.22] & [60.18] & [61.92] & [64.44] & [65.38] & [43.90] \\ \hline 
$n=200$, $a=1.5$ in (\ref{eg52}) & 100 & 100 & 100 & 100 & 100 & 100 & 100 \\ 
 & [100] & [100] & [100] & [100] & [100] & [100] & [100] \\ \hline \hline
 \end{tabular}
\end{center}
 \end{table}

\section{Application to German Bundesbank Data}

Banking throughout the world is based on credit, or on trust in the debtor's ability to fulfill his/her debt obligation. However, facing increasing pressure from markets and regulators, banks have based their risk analysis, increasingly, on statistical techniques to judge or predict corporate bankruptcy. This is known as rating or scoring. Its main purpose is to estimate the financial status of a company and, if possible, to estimate the probability of a company default on its debt obligations within a certain period. Logistic regression is probably the most commonly used technique to model the probability of default  and logistic partially linear models may also be advantageous because of its flexibility, in allowing for the possibly nonlinear effects of one continuous covariate.

We apply the methodology  to the   German Bundesbank Data in year 2002. 
The data provided by CRC 649, Humboldt-Universit\"{a}t zu Berlin,   contained  6123 companies  of which 186 were insolvent. Each firm is described by 28 financial ratio variables, $x1, \dots, x28$,  and those of insolvent firms were collected two years prior to insolvency.   
 To ensure the value of some variables as the denominator should not be zero when calculating the ratios, 2079 firms were retained with 92 insolvent. 
Though removing almost two thirds of the sample may seem excessive,
we did not intend to analyze the majority of  firms in the database.  
 The focus was to investigate (i) differences between the financial ratios of the
solvent and insolvent firms, and (ii) how  the nonlinear effects  improve parametric logistic fitting.

Based on support vector machines and for a much larger data sample spanning from 1996 through to 2002, Chen et al. (2011)   selected $x24$ (accounts payable/sales) 
 measuring account payable turnover, as the best predictor, and subsequently selected  $x3$ (operating income/total assets)
measuring profitability, $x15$ ((cash and cash equivalents)/total assets)
measuring liquidity, $x12$ (total liabilities/total assets)
measuring leverage, $x26$ (increase (decrease) inventories/inventories)
measuring percentage of incremental inventories, $x22$ (inventories/sales)
measuring inventory turnover, $x5$ ((earnings before interest and tax)/total assets)
and $x2$ (net income/sales)
measuring net profit margin. For year 2002 data, we found that $x3$ and $x5$ have a large sample correlation coefficient  0.95 and thus $x5$ is removed from our analysis  and we further  include $x25$ (log(total assets)) measuring firm size, as it is shown to be an important variable on predicting the probability of bankruptcy in the literature (see, e.g., Lopez 2004).  In summary, there are 8 predictors, $x2$, $x3$, $x12$, $x15$, $x22$, $x24$, $x25$, and $x26$, and  a binary response. See  Chen et al. (2011) for detail descriptions about the data.

Since $x24$ was selected as the most important predictor by Chen et al. (2011), we model its effects nonparametrically, while retaining linear trends for the remaining predictors in a logit model. 
The variable $x24$ measuring account payable turnover is a short-term liquidity measure for quantifying the rate at which a firm pays off its suppliers. Generally speaking, ``the firms with higher account payable turnover  will have less ability to convert their accounts into sales, have lower revenues, and go bankrupt more readily" (Chen et al. 2011). However this measure is specific to different industries; every industry has a slightly different standard. 
Further examination of $x24$ indicates that most values lie in $[0, 0.5]$ with only 15 observations in $(0.5, 20.52)$. 
 If  those 15 observations are excluded, then the  sample size becomes 2064, in which 91 are insolvent.  An alternative approach, suggested by a reviewer, is taking logarithm of $(x24+0.001)$ (0.001 is added since $x24$ includes 0's) 
 and retaining the sample size $n=2079$.

Local linear smoothing with  the Epanechnikov kernel  is used. The values of bandwidth for $x24$,  0.125,  0.1, and 0.08, are equally-spaced on a  logarithmic  scale,  corresponding to df  4.94, 5.92, and 7.17 respectively. The bandwidth that minimizes AICc is 0.125 and $h_{HS}=0.1$.
The curves for $m(x24)$ with pointwise confidence intervals based on empirical Fisher information matrices  are shown in Figure 2 with $h=0.125$ and 0.1,  and 
the proposed test for testing $H_0: m(x24)$ is a constant, gives highly significant $p$-values $< 10^{-14}$ for all 3 values of the bandwidth, indicating significance of $m(x24)$ in predicting probability of default .
A linear logistic model gives a positive slope 10.37 for $x24$ with a highly significant $p$-value $< 10^{-15}$. Since Chen et al. (2011) interpreted the linear trend as  higher default probability with  high turnover, we attempt to interpret the seemingly non-linear curves in Figure 2 as follows.
Taking the curve with $h=0.1$ in Figure 2, when $x24$ increases from 0.1 to 0.3, the  estimate increases about 1.895, which means the odds ratio for a firm with $x24=0.3$ to become insolvent  is $\exp(1.895)=6.653$ times relative to that for a firm with $x24=0.1$. On the other hand,  between $x24=0.3$ and $0.4$, the  estimate decreases by an amount of $-0.607$, implying that the odds ratio for a firm with $x24=0.4$ to become insolvent  is $\exp(-0.607)=0.545$ times relative to that for a firm with $x24=0.3$. Thus our analysis gains new insight  suggesting  that a German firm is likely to go bankrupt when it has a higher turnover for roughly 97.5\% of firms (0.3 is approximately 97.5-percentile of $x24$), but for those firms with $0.3<x24<0.4$ (approximately 97.5- to 99-percentile),  the default probability decreases as $x24$ increases. 

\begin{figure}
\includegraphics[width=6.5in,height=3.5in]{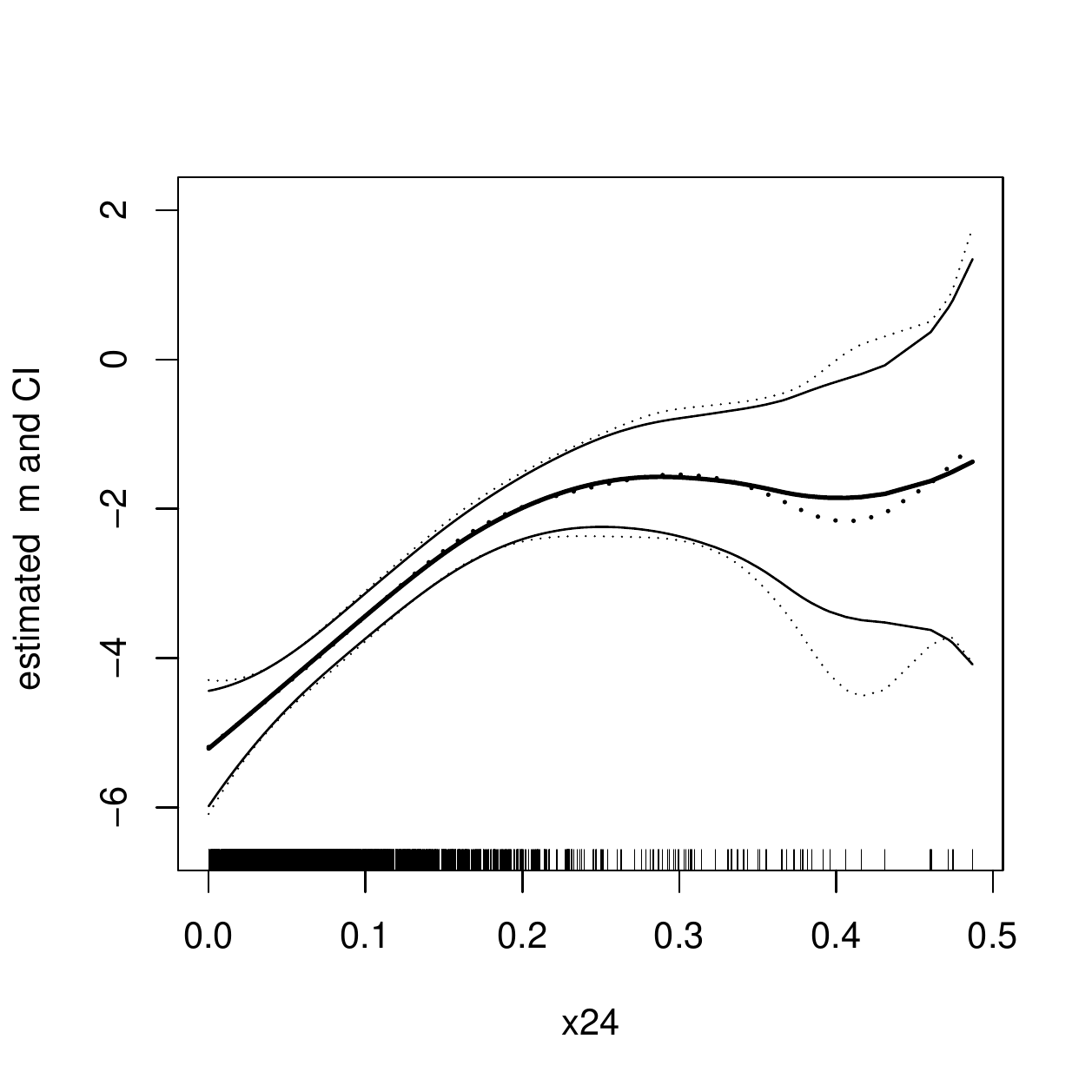} 
\caption{Plot of the nonlinear trends of $x24$ in predicting the probability of bankruptcy using bandwidth $h=0.125$ (solid line) and $h=0.1$ (dash line) with 95\% pointwise confidence intervals for the 2002 German Bundesbank Data.}
\end{figure}

\begin{figure}
\includegraphics[width=6.5in,height=3.5in]{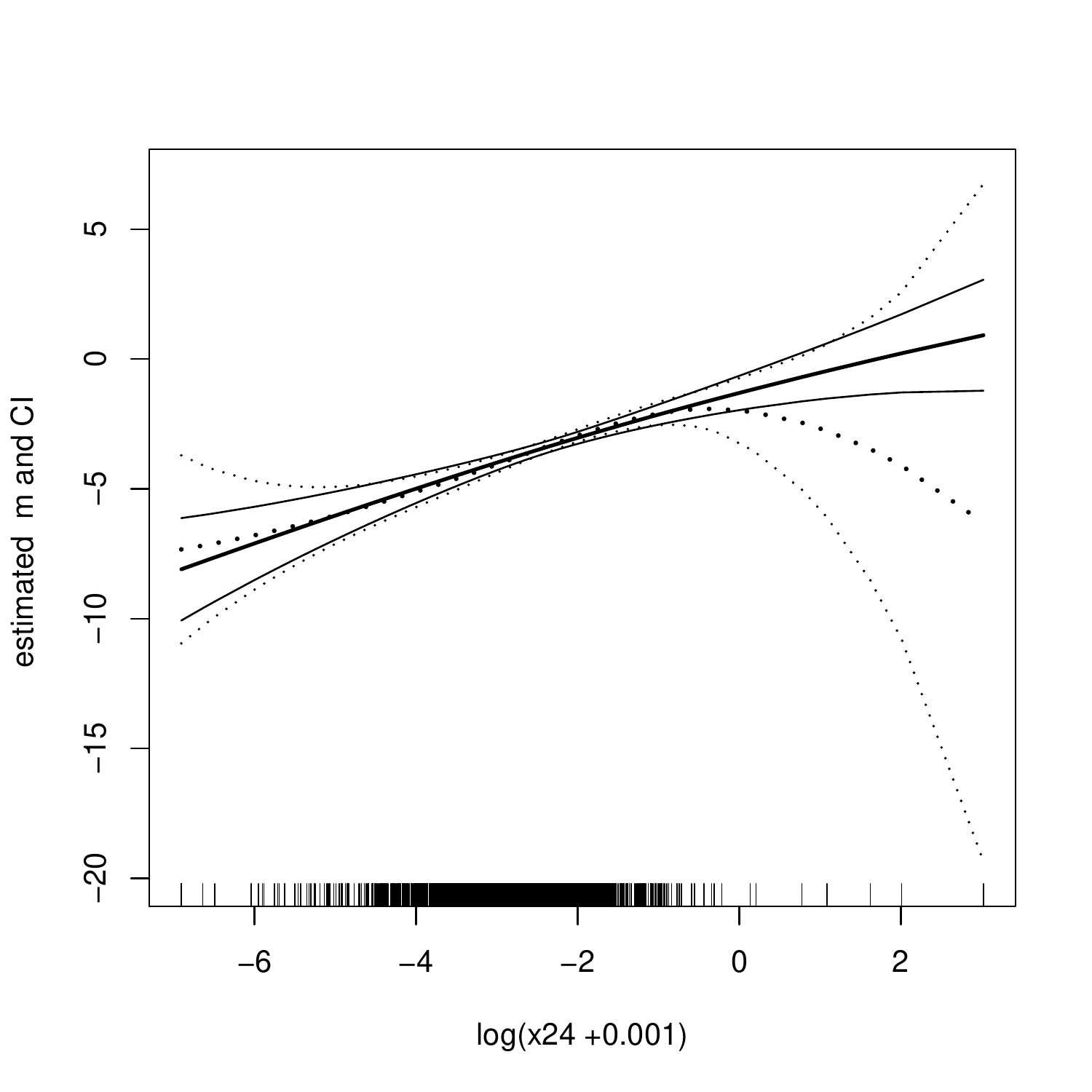} 
\caption{Plot of the  trends of $logx24$ in predicting the probability of bankruptcy using bandwidth $h=6$ (solid line) and $h=3.375$ (dash line) with 95\% pointwise confidence intervals for the 2002 German Bundesbank Data.}
\end{figure}

For smoothing on log$(x24+0.001)=logx24$ with $n=2079$, the values of $h$ are  $6,$ 4.5, and 3.375,   corresponding to df   2.32,  2.83, and  3.69  respectively.  The bandwidth that minimizes AICc is  3.375 and $h_{HS}=6$.
The curves for $logx24$  with $h=3.375$ and 6 shown in Figure 3 have a  linear tendency for $logx24 <-0.2$ $(x24< 0.8)$, and the confidence intervals corresponding to $h=3.375$ imply some uncertainty near the right-hand end points.
The tests for $H_0: m(logx24)$ is a constant, give highly significant $p$-values of $< 1.5 \times 10^{-12}$. If using a linear trend for $logx24$, the slope is 0.927 with a significant $p$-value of $2.0 \times10^{-12}$.  The analysis using $logx24$ implies that a German firm is likely to go bankrupt when it has a high turnover in the log scale, but the linear trend is uncertain for those with $x24>0.8$. Hence  an interpretation in Chen et al. (2011)  that a German firm is likely to go bankrupt when it has high turnover may not be entirely correct;    the  effects of $x24$ on the probability of bankruptcy may be nonlinear for those with large turnovers, as shown in Figures 2 and 3.  

\section{Discussion}

We develop local and global analysis of deviance expressions and associated integrated likelihood ratio tests for generalized partially linear models with canonical links based on fitting local $p$-th order polynomials. 
Though the idea of nonparametric analysis of deviance  is not
new (Hastie and Tibshirani, 1990), the work in this paper provides theoretical justifications  that connect to the classic framework.
Theorems 2 and 4(b) are restricted for $p=1$ and 3 only, while Theorems 1, 3, 4(a)(c), and 5 are for a nonnegative integer $p$. 
As a  by-product,  new estimators for the canonical parameter and response mean  are proposed   and Theorems 2 and 4(b)  show that the integrated likelihood quantities are asymptotically global likelihood quantities with the new estimators. The new estimator  $\hat{\theta}_i^*$ or 
$\breve{\theta}_i^*$ for the canonical parameter is formed by
combining  locally fitted $\hat{\theta}_i(x)$ or $\breve{\theta}_i(x)$  through weighted integration and thus utilize all locally fitted parameters, which is different from the conventional approach of focusing on $\hat{\beta}_0$. The integrated likelihood approach of combining local likelihood appears to be new in the smoothing literature, though it was discussed by   Severini (2007) and Lehmann (2006) in different settings. 
The numerical results of $n=100$ and $200$  show that the test statistics under the null hypothesis follow the asymptotic $\chi^2$-distribution reasonably well  and the performance under alternative hypotheses is sometimes more powerful than Wood (2013)   in  the R package {\tt mgcv}. 
It has been suggested by a reviewer to investigate the asymptotic power of the proposed tests.
Since there is no simple explicit expression for Fisher
information for generalized linear models (\ref{glm}), we conjecture that the study of power may be focused on special cases of 
logistic and Poisson models, which will be explored for future research.
For a smaller sample size such as $n=50$, two reviewers has suggested to develop a bootstrap procedure for calculating the sample critical values for further investigation.

The local analysis of deviance  in  (\ref{laod}) and  in Theorem 4(a) are derived assuming a fixed value of bandwidth. It is straightforward to obtain local analysis of deviance expressions with varying values of bandwidth at different $x$, but how to combine them to form global analysis of deviance  will be an interesting problem. 
Like all smoothing-based tests, the $p$-values of the integrated likelihood ratio tests depend on the values of the smoothing parameter. We recommend plotting the ``significant trace" (Bowman and Azzalini 1998) to assess the evidence across a wide range of values of $h$ and looking for some overall trends.
 For fitting generalized partially linear models, a practical problem is how to choose the predictor to be modelled nonparametrically. One approach may be  based on selecting the most significant predictor based on the smallest $p$-value of integrated likelihood ratio tests when using approximately the same degrees of freedom for smoothing. This idea and the related model selection problems with a diverging  number of linear covariates (Wang et al. 2014) may be explored for future research.   A topic for further investigation is the problem of   bandwidth selection for models (\ref{nonpar}) and (\ref{model})  based on  cross-validating the deviance or minimizing the Akaike information
criterion.
Further extension on  developing analysis of deviance  for generalized partially linear models with non-canonical links, for multiplicative bias reduction methods (Kosmidis and Firth 2009),  for hazard estimation (Nielsen and Tanggaard 2001) as the proportional hazards models and Poisson regression are connected,  and for generalized additive models with multiple nonparametric functions remain to be investigated.

\noindent {\bf Acknowledgement}\\
We thank the editor, associate editor , and two anonymous referees for their constructive comments and suggestions.
 The work was conceived during the visit of the second author to the  Humboldt-Universit\"{a}t zu Berlin supported
by  CRC 649 ``Economic Risk."  The support is greatly appreciated.
The first author was partially supported by SKBI School of Business, Singapore Management University,  
the second author (corresponding author) was partially supported by  the Ministry of Science and Technology NSC 101-2118-M-007-002-MY2 in Taiwan, and   both
  authors were partially supported by    CRC 649  ``Economic Risk."  The authors wish to thank Mr. Leslie  Udvarhelyi who assisted in the proof-reading of the manuscript.

\appendix
\section*{Appendix}

\noindent  {\sc Conditions (A).}
\begin{enumerate}

\item[(A1).] The kernel $K(\cdot)$ is a Lipschitz continuous, bounded and
            symmetric probability
            density function, having a support on a compact interval, say $[-1,1]$. 
For  the $x_i$s that fall in the conventional boundary region, a boundary-corrected kernel
is used to ensure  $\int K_h(x_i-x) dx=1$.

\item[(A2).]  The covariate $X$ is assumed to have a bounded support  $\mathcal{X}$.

\item[(A3).] The function $(\partial^2/\partial \mu^2) \ell \{ y; G(\mu)\} <0$ for $\mu \in \mathcal{R}$ and $y$ in the range of the response variable.

\item[(A4).] The functions $L^{\prime}$, $\theta^{(p+2)}$, $b^{\prime \prime}(\theta(\cdot)) \equiv V(\cdot)$, $V^{\prime \prime}$, and $G^{(3)}$ are continuous with respect to $x$.

\item[(A5).]  For each $x$ in $\mathcal{X}$, 
$V(x)$ and $G^{\prime}(\mu(x))$ are non-zero.

\item[(A6).] For  (\ref{model}), the covariate vector $\bz$ is assumed to have a bounded support.

\end{enumerate}

Conditions  (A1)-(A5) are similar to those in Fan et al. (1995). Without loss of generality, Condition (A2)   is  satisfied in practice by  strictly increasing transformations of $X$ when the support of $X$ before transformation is unbounded. Similar explanations can be said about condition (A6). Condition (A3)
ensures that the local polynomial estimate $\hat{\beta}$ lies in a compact set. Conditions (A3) and (A5) imply that $V(x)>0$ for $x \in \mathcal{X}$. Condition (A4) implies that all thrid derivatives of $\ell \{ y; G(\mu)\}$ with respect to $\mu$  are continuous, and $V^{\prime}$ and $\mu^{\prime}$ are continuous with respect to $x$.

\noindent  {\bf Proof of Proposition}

\noindent  When $p=1$,  the bias of $\theta_i^*$ is expressed as follows:
\begin{eqnarray}  
\nonumber \lefteqn{  \E \left \{ \int (\hat{\beta}_0(x) + (x_i-x) \hat{\beta}_1(x)) K_h(x_i-x) dx \right\} - \theta_i =} \\
\nonumber & \E \left [  \int \{ (\hat{\beta}_0(x)- \beta_0(x)) + (x_i-x) (\hat{\beta}_1(x) - \beta_1(x))\} K_h(x_i-x) dx - \right.\\
 & \left.  \int \{\beta_2(x)   (x_i-x)^2 + r(x, x_i)\}  K_h(x_i-x) dx \right ],
\label{p1}
\end{eqnarray}
where $r(x, x_i)$ denotes the remainder terms.   Plugging the first-order term of  the asymptotic bias of  $\hat{\beta}_0(x)$  (Fan et al. (1995) and Aerts and Claeskens (1997)) in  (\ref{p1}), leads to cancellation with the $\beta_2(x)$-term  in  (\ref{p1}). The  remaining  term $\int (x_i-x) (\hat{\beta}_1(x) - \beta_1(x)) K_h(x_i-x) dx$ is of order  $h^4$. Thus   the $h^2$-order term in (\ref{p1}) is   zero. Similar arguments can be shown for $p=3$.

Fan et al. (1995) and Aerts and Claeskens (1997) showed that the variance of $\hat{\beta}_j(x)$ is of order $n^{-1}h^{-2j-1}$ when $p$ is odd.  Then the variance of $\{\hat{\beta}_0(x) + (x_i-x) \hat{\beta}_1(x)\}$ is of order $n^{-1} h^{-1}$ and hence the variance of $\theta_i^*$ is of order $n^{-1} h^{-1}$. Finally, it is straightforward to show (c)  based on (a)  and (b) since $\mu_i^*=G^{-1}(\theta_i^*)$.

\noindent  {\bf Proof of Theorem 2}

\noindent
We only need to show Theorem 2(a) while Theorem 2(b) follows directly from Theorem 2(a). 
For the left-hand side  of  (\ref{eq2}), ignoring the $c(y, \phi)$ and $a(\phi)$ terms in (\ref{glm}) which is unrelated to $x$,  the integrated likelihood is
$$\int  \ell_x(\by; \hat{\btheta}_x) dx=  \sum_i \left\{ y_i \int \hat{\theta}_i (x)K_h (x_i -x)  dx -\int b(\hat{\theta}_i (x))  K_h (x_i -x) dx \right\}.$$
By  Taylor's expansion,
\[  b(\hat{\theta}_i (x)) = b(\theta_i^*) + b^{\prime} (\theta_i^*) \{ \hat{\theta}_i (x)- \theta_i^* \} +  b^{\prime \prime} (\theta_i^*) \{  \hat{\theta}_i (x)- \theta_i^* \}^2 /2+ r_i(x),\]
where $r_i(x)$ denotes the remainder terms.
For the linear term,
\[ \int b^{\prime} (\theta_i^*) ( \hat{\theta}_i (x)- \theta_i^* )K_h (x_i -x)  dx
=b^{\prime} (\theta_i^*) ( \theta_i^* -  \theta_i^*) =0.\]
The quadratic term 
$\int \{ \hat{\theta}_i (x)- \theta_i^* \}^2K_h (x_i -x)  dx$
$=  \int \hat{\theta}_i (x)\{\hat{\theta}_i (x)-\theta_i^{*}\}   K_h (x_i -x)  dx.$ For $\hat{\theta}_i (x)-\theta_i^{*}=
(\hat{\theta}_i (x)- \theta_i) - (\theta_i^{*}-\theta_i)$,
the first term $\hat{\theta}_i (x)- \theta_i= (\hat{\beta}_0 - \beta_0) + (\hat{\beta}_1 - \beta_1)(x_i -x) + \dots + (\hat{\beta}_p - \beta_p) (x_i -x)^p +
r^{\prime}_i(x)$  by  Taylor's expansion of $\theta_i$, where $r^{\prime}_i(x)$ denotes the remainder terms. Then based on Theorem 1(a) of Fan et al. (1995), 
$ \hat{\theta}_i (x)- \theta_i$ is of order  $\CO(h^{p+1})$. For the second term
$ \theta_i^{*}-\theta_i$, it is of order $\Co(h^{(p+1)})$ based on the Proposition. 
Thus $ \int b(\hat{\theta}_i (x))  K_h (x_i -x) dx = b(\theta_i^*) + \CO(h^{(p+1)})$ and   (\ref{eq2})  is proved.

\noindent  {\bf Proof of Theorem 3}

\noindent
Let $\ell(y_i; a_0)$ be the likelihood corresponding to $y_i$ with  $\theta_i= a_0$.
Define $\ell_x(\by; a_0)= \sum_i \ell(y_i; a_0) K_h(x_i-x)$,  and  it is clear that $\int \ell_x(\by; a_0) dx= \ell(\by; a_0)$ under (A1).  
Recall that $\hat{\bbeta}=(\hat{\beta}_0, \dots, \hat{\beta}_p)^{\top}$ maximizes local likelihood at $x$, $\ell_x (\by; \btheta_x)$, with local polynomial approximation. Expanding $\ell_x (\by; \hat{\btheta}_x)$, which is a function of $\hat{\bbeta}$,  around a $(p+1)$-length vector $\ba_0=(a_0, 0, \dots, 0)^{\top}$,
\begin{equation}
\ell_x (\by; \hat{\btheta}_x) - \ell_x(\by; \ba_0)=\left  \{ \frac{\partial \ell_x}{\partial \bbeta} (\by; \ba_0 )\right \}^{\top} ( \hat{\bbeta} - \ba_0 ) + \frac{1}{2}  ( \hat{\bbeta} - \ba_0 )^{\top}  \frac{\partial^2 \ell_x}{\partial \bbeta^2} (\by; \ba_0) ( \hat{\bbeta} - \ba_0 )
 + \CO_p(\sqrt{n^{-1} h^{-1}}),
\label{eqt3.1}
\end{equation}
where $\partial \ell_x(\by; \ba_0 )/ \partial \bbeta=  \partial \ell_x(\by; \btheta_x )/ \partial \bbeta \bigm|_{\ba_0}$ and similarly for
$\partial^2 \ell_x (\by; \ba_0)/ \partial \bbeta^2$.
Substituting the expansion
\[ \hat{\bbeta} - \ba_0 = i_x(\ba_0)^{-1} \frac{\partial \ell_x}{\partial \bbeta} (\by; \ba_0)+ \CO_p(\sqrt{n^{-1} h^{-1}}),\]
where $i_x(\ba_0) =\E \left \{ - \partial^2 \ell_x(\by; \btheta_x)  /\partial \bbeta^2 \right \} (\ba_0)$,
we have for (\ref{eqt3.1}),
\begin{eqnarray*}
 \left  \{ \frac{\partial \ell_x}{\partial \bbeta} (\by; \ba_0) \right \}^{\top}   i_x(\ba_0)^{-1} \frac{\partial \ell_x}{\partial \bbeta} (\by; \ba_0) + \frac{1}{2} \left\{  \frac{\partial \ell_x}{\partial \bbeta} (\by; \ba_0) \right \}^{\top}   i_x(\ba_0)^{-1}   \frac{\partial^2 \ell_x}{\partial \bbeta^2} (\by; \ba_0)   i_x(\ba_0)^{-1} \frac{\partial \ell_x}{\partial \bbeta} (\by; \ba_0)  \\
 + \CO_p(\sqrt{n^{-1} h^{-1}}).
\end{eqnarray*}
Since $\partial^2 \ell_x(\by; \ba_0) / \partial \bbeta^2 = - i_x(\ba_0)  + \CO_p(\sqrt{nh})$, 
it follows that $\int 2  \{ \ell_x (\by; \hat{\btheta}_x) - \ell_x(\by; a_0) \} dx$ is
\begin{eqnarray*}
  \int \left \{ \frac{\partial \ell_x}{\partial \bbeta} (\by; \ba_0) \right \}^{\top}   i_x(\ba_0)^{-1} \frac{\partial \ell_x}{\partial \bbeta} (\by; \ba_0)  dx 
=   \left \{i(a_0)^{-1/2} \frac{\partial \ell}{\partial \btheta} (\by; a_0) \right \}^{\top}  H_p^* \left\{  i(a_0)^{-1/2} \frac{\partial \ell}{\partial \btheta} (\by; a_0) \right\},
\label{eq3.2}
\end{eqnarray*}
where the last expression is obtained by plugging in  the explicit expressions of $ \partial \ell_x(\by; \ba_0)/ \partial \bbeta $ and $ i_x(\ba_0)^{-1}$,  and $i(a_0)$  is the information matrix under $H_0$. 

By standard ML theory,   $2\{ \ell(\by; \hat{a}_0)- \ell(\by; a_0) \}$ is asymptotically distributed according to a $\chi^2$ distribution with $1$ degree of freedom and $i(a_0)^{-1/2}  \partial \ell(\by; a_0) / \partial \btheta $ is asymptotically normally distributed with mean vector 0 and identity covariance matrix.  Then the test statistic (\ref{liktest}) becomes 
\[   \left \{i(a_0)^{-1/2} \frac{\partial \ell}{\partial \btheta} (\by; a_0) \right \}^{\top}  \{ H_p^*-  P_0 \} 
\left\{  i(a_0)^{-1/2} \frac{\partial \ell}{\partial \btheta} (\by; a_0) \right\},\]
 where   $P_0$ is an $n\times n$ matrix with $1/n$ in all entries. 
From Huang and Chen (2008), 
 $H_p^*- P_0$ is symmetric and asymptotically idempotent. Thus  the test statistic (\ref{liktest}) has an asymptotic $\chi^2$-distribution with df $\text{tr}(H_p^*) -1$.

\noindent  {\bf Proof of Theorem 5}

\noindent
The proof is an extension from  that  of Theorem 3. Let $\balpha_0= (a_0, \alpha^{\top})^{\top}$ denote the parameter vector under $H_0$ and
 $\ell(\by;  \balpha_0)$ be the corresponding likelihood.
Define local likelihood at $x$ under $H_0$ by  $\ell_x(\by; \balpha_0)= \sum_i \ell(y_i; \balpha_0) K_h(x_i-x)$ and  hence $\int \ell_x(\by; \balpha_0) dx= \ell(\by; \balpha_0)$.  
We consider
$\int  \ell_x(\by; \breve{\btheta}_x) dx - \ell (\by; \balpha_0) $ and  
$ \ell (\by; \hat{\balpha}_0)    - \ell(\by; \balpha_0)$
separately, whose difference becomes (\ref{liktest}).

 Let $\bb=( \beta^{\top}, \alpha^{\top})^{\top}$  be the parameter vector under $H_1$.  
Expanding $\ell_x (\by; \breve{\btheta}_x)$, which is a function of $\breve{\bb} =(\breve{\bbeta}^{\top}, \breve{\alpha}^{\top})^{\top}$,  around a $
(p+K+1)$-length vector $\ba_0=(a_0, 0, \dots, 0, \alpha^{\top})^{\top}$,
\begin{equation}
\ell_x (\by; \breve{\btheta}_x) - \ell_x(\by; \ba_0)=\left  \{ \frac{\partial \ell_x}{\partial \bb} (\by; \ba_0) \right \}^{\top} ( \breve{\bb} - \ba_0 ) + \frac{1}{2}  ( \breve{\bb} - \ba_0 )^{\top}  \frac{\partial^2 \ell_x}{\partial \bb^2} (\by; \ba_0) ( \breve{\bb}- \ba_0 )
 + \CO_p(\sqrt{n^{-1} h^{-1}}).
\label{eqt5.1}
\end{equation}
Substituting the expansion under $H_0$,
\[ \breve{\bbeta} - \ba_0 = i_x(\ba_0)^{-1} \frac{\partial \ell_x}{\partial \bb} (\by; \ba_0)+ \CO_p(\sqrt{n^{-1} h^{-1}}),\]
where $i_x(\ba_0) =\E \{ - \partial^2 \ell_x(\by; \btheta)/  \partial \bb^2   \} (\ba_0)$,
we have for (\ref{eqt5.1}),
\begin{eqnarray*}
 \left  \{ \frac{\partial \ell_x}{\partial \bb} (\by; \ba_0) \right \}^{\top}   i_x(\ba_0)^{-1} \frac{\partial \ell_x}{\partial \bb} (\by; \ba_0) + \frac{1}{2} \left\{  \frac{\partial \ell_x}{\partial \bb} (\by; \ba_0) \right \}^{\top}   i_x(\ba_0)^{-1}   \frac{\partial^2 \ell_x}{\partial \bb^2} (\by; \ba_0)   i_x(\ba_0)^{-1} \frac{\partial \ell_x}{\partial \bb} (\by; \ba_0)  \\
 + \CO_p(\sqrt{n^{-1} h^{-1}}).
\end{eqnarray*}
Since $\partial^2 \ell_x(\by; \ba_0) /\partial \bb^2 = - i_x(\ba_0)  + \CO_p(\sqrt{nh})$ under $H_0$, 
it follows that $\int 2  \{ \ell_x (\by; \breve{\btheta}_x) - \ell_x(\by; \balpha_0) \} dx$ is asymptotically
\begin{eqnarray*}
       \left \{i(\balpha_0)^{-1/2} \frac{\partial \ell}{\partial \btheta} (\by; \balpha_0) \right \}^{\top}  \{ H_p^* + P_z\}  \left\{  i(\balpha_0)^{-1/2} \frac{\partial \ell}{\partial \btheta} (\by; \balpha_0) \right\},
\end{eqnarray*}
where $P_z$ is the projection matrix for $Z$  and $i(\balpha_0)$  is the information matrix under $H_0$.  The last expression  is obtained by plugging in  the explicit expressions of 
$ \partial \ell_x(\by; \ba_0)/ \partial \bb $ and $ i_x(\ba_0)^{-1}$ and using the assumption that $\bx$ and $Z$ are orthogonal.

 For the other term, $2\{ \ell(\by; \hat{\balpha_0})- \ell(\by; \balpha_0) \}$ is asymptotically distributed according to a $\chi^2$-distribution with $(K+1)$ degree of freedom under $H_0$ by standard ML theory.  Hence the test statistic (\ref{liktest}) becomes 
\[   \left \{i(\balpha_0)^{-1/2} \frac{\partial \ell}{\partial \btheta} (\by; \balpha_0) \right \}^{\top}  \{ H_p^*-  P_0 \} 
\left\{  i(\balpha_0)^{-1/2} \frac{\partial \ell}{\partial \btheta} (\by; \balpha_0) \right\},\]
 where   $P_0$ is the same as in the proof of Theorem 3.  Also 
$  i(\balpha_0)^{-1/2}  \partial \ell (\by; \balpha_0)/ \partial \btheta  $ is asymptotically normally distributed with mean vector 0 and identity covariance matrix. Hence the test statistic (\ref{liktest}) has an asymptotic chi-square distribution with df $\text{tr}(H_p^*) -1$.

\subsection*{References}
\setstretch{1}

\begin{description}

\item Aerts, M., and Claeskens, G. (1997), ``Local Polynomial Estimation in Multiparameter Likelihood Models,"
    {\it Journal of the American Statistical Association},  92, 1536-1545.

\item  Bowman, A.W., and Azzalini, A. (1997),
 {\it  Applied Smoothing Techniques for Data Analysis,}
  London: Oxford.

\item Carroll, R. J., Fan, J., Gijbels, I., and Wand, M. P. (1997),
``Generalized Partially Linear Single-index Models," {\it Journal of the American Statistical Association},  92, 477-489.

\item Chen, S., H\"ardle, W.,  and Moro, R. (2011),
``Modeling Default Risk with Support Vector Machines,"
{\it Quantitative  Finance,}  11,   135-154.

 \item Fan, J., and Gijbels, I. (1996), {\it Local Polynomial 
Modelling and Its Applications,}
London:  Chapman and Hall.

\item  Fan, J.,  Heckman, N. E., and  Wand, M. P. (1995), ``Local Polynomial Kernel Regression for Generalized Linear Models and Quasi-Likelihood Functions," {\it Journal of the American Statistical Association},  
 90, 141-150.

\item Fan, J., Zhang, C., and Zhang, J. (2001), ``Generalized Likelihood Ratio Statistics and Wilks Phenomenon,"
 {\it Annals of Statistics}, 29, 153-193.

\item Green, P. J., and Silverman, B. W. (1994), {\it 
Nonparametric Regression and Generalized Linear Models: a Roughness Penalty 
Approach,} London: Chapman  and  Hall.

\item  H\"ardle, W., Muller, M.,  and Mammen, E.  (1998),
``Testing Parametric Versus Semiparametric Modeling in Generalized Linear Models,''
{\it Journal of the American Statistical Association}, 93, 1461-1474.

\item   H\"ardle, W. K., M\"uller, M., Sperlich, S., and Werwatz, A. (2004),
  {\it Nonparametric and Semiparametric Models}, Berlin: Springer.

\item Hastie, T. J., and Tibshirani, R. J. (1990), {\it Generalized Additive 
Models,} London:
Chapman and Hall.

\item Hastie, T. J., and Tibshirani, R. J. (1987), ``Local Likelihood Estimation,"
{\it Journal of the American Statistical Association},   82, 559-567.

\item Horowitz, J. L., and Spokoiny, V. G. (2001), ``An Adaptive, Rate--optimal Test of a Parametric Mean‐-Regression Model Against a Nonparametric Alternative," {\it Econometrica}, 69, 599-631.

\item  Huang, L.-S., and Chan, K.-S. (2014),
``Local Polynomial and Penalized Trigonometric Series
Regression," {\it Statistica Sinica,}  24, 1215-1238.

\item Huang, L.-S.,  and Chen, J. (2008),  ``Analysis of Variance,
Coefficient of Determination, and F-test for Local Polynomial Regression,''
{\it Annals of Statistics}, 36, 2085-2109.

\item  Huang, L.-S., and Davidson, P. W.  (2010),
``Analysis of Variance and  $F$-tests for Partial Linear Models with Applications to  Environmental Health Data,"
{\it Journal of the American Statistical Association}, 
 105,  991-1004.

\item Hurvich, C.M., Simonoff, J.S., and  Tsai, C.-L.  (1998),
``Smoothing Parameter Selection in Nonparametric
  Regression Using an Improved Akaike Information Criterion,''
 {\it Journal of the Royal Statistical Society}, Series B,    60, 271-293.

\item Kosmidis, I., and Firth, D. (2009), ``Bias Reduction in Exponential Family Nonlinear Models," {\it Biometrika},  96, 793-804.

\item Lehmann, E. L. (2006), ``On  Likelihood Ratio Tests," in {\it Optimality: The Second Erich L. Lehmann Symposium,}
Institute of Mathematical Statistics Lecture Notes - Monograph Series  vol. 49, ed J. Rojo,   Beachwood, OH: Institute of Mathematical Statistics, pp 1-8.

\item Li, R., and Liang, H. (2008), ``Variable Selection in Semiparametric Regression Modeling,"  
{\it Annals of Statistics}, 36, 261-286.

\item Loader, C. (1999), {\it Local Regression and Likelihood,}  New York: Springer.

\item Lopez, J. A. (2004), ``The Empirical Relationship Between Average Asset Correlation, Firm Probability of Default, and Asset Size,"
 {\it Journal of Financial Intermediation}, 13, 265-283.

\item Nielsen, J. P., and Tanggaard, C. (2001), ``Boundary and Bias Correction in Kernel Hazard Estimation,"
 {\it Scandinavian Journal of Statistics,} 28, 675-698.

\item McCullagh, P., and   Nelder, J. A. (1989), {\it Generalized Linear Models,} (2nd ed.) London: Chapman and Hall.

\item Severini, T. A. (2007), ``Integrated Likelihood Functions for Non-Bayesian Inference," {\it Biometrika,} 94, 529-542.

\item Simon, G. (1973), ``Additivity of Information in Exponential Family Probability Laws,"
{\it Journal of the American Statistical Association},   68, 478-482.

\item Wang, L., Xue, L., Qu, A., and Liang, H. (2014), ``Estimation and Model Selection in Generalized Additive Partial Linear Models for Correlated Data with Diverging Number of Covariates," 
{\it Annals of Statistics}, 42, 592-624.

\item Wood, S. N. (2006), {\it Generalized Additive Models: An Introduction with R}, Boca Raton, FL: Chapman and Hall/CRC Press. 

\item Wood, S. N. (2013), ``On $p$-values for Smooth Components of an Extended Generalized Additive Model," 
{\it Biometrika,}  100, 221-228.

\end{description}

\end{document}

%% file: defs.tex
\newcommand{\CO}{{\mathcal{O}}}

\newcommand{\E}{\mbox{\sf E}}     

\def\Co{{\scriptstyle \mathcal{O}}} 
\def\CO{\mathcal{O}}
\def\Co{{\scriptstyle{\mathcal{O}}}}

